\documentclass[
11pt,                          
english                        
]{article}

\usepackage[english]{babel}    
\usepackage{amsmath}           
\usepackage{amssymb}           
\usepackage[utf8]{inputenc}    
\usepackage[T1]{fontenc}       
\usepackage{eucal}             
\usepackage{bbm}               
\usepackage{exscale}           
\usepackage{mathtools}         
\usepackage[amsmath,thmmarks,hyperref]{ntheorem} 
\usepackage[sort]{cite}        
\usepackage{array}             
\usepackage{stmaryrd}          
\usepackage{paralist}          
\usepackage[a4paper]{geometry} 
\usepackage{xspace}            
\usepackage{ifdraft}           
\usepackage[expansion=false]{microtype}          
\usepackage[nottoc]{tocbibind} 
\usepackage[backref=page,final=true]{hyperref}   

\ifdraft{\synctex=1}{}

\renewcommand{\mathbb}[1]{\mathbbm{#1}}


\newcommand{\refitem}[1] {\textit{\ref{#1}.)}}

\numberwithin{equation}{section}

\allowdisplaybreaks

\let\originalleft\left
\let\originalright\right
\renewcommand{\left}{\mathopen{}\mathclose\bgroup\originalleft}
\renewcommand{\right}{\aftergroup\egroup\originalright}

\theoremheaderfont{\normalfont\bfseries}
\theorembodyfont{\itshape}
\newtheorem{lemma}{Lemma}[section]
\newtheorem{proposition}[lemma]{Proposition}
\newtheorem{theorem}[lemma]{Theorem}
\newtheorem{corollary}[lemma]{Corollary}
\newtheorem{definition}[lemma]{Definition}

\theorembodyfont{\rmfamily}

\newtheorem{example}[lemma]{Example}
\newtheorem{remark}[lemma]{Remark}

\makeatletter
\def\theorem@checkbold{}
\makeatother

\theoremheaderfont{\scshape}
\theorembodyfont{\normalfont}
\theoremstyle{nonumberplain}
\theoremseparator{:}
\theoremsymbol{\hbox{$\boxempty$}}
\newtheorem{proof}{Proof}

\newenvironment{examplelist}{\begin{compactenum}[\itshape i.)]}{\end{compactenum}}
\newenvironment{lemmalist}{\begin{compactenum}[\itshape i.)]}{\end{compactenum}}
\newenvironment{theoremlist}{\begin{compactenum}[\itshape i.)]}{\end{compactenum}}
\newenvironment{propositionlist}{\begin{compactenum}[\itshape i.)]}{\end{compactenum}}

\newcommand{\cc}[1]          {\overline{{#1}}}
\newcommand{\at}[1]          {\big|_{#1}}
\newcommand{\argument}       {\,\cdot\,}
\newcommand{\earlier}        {\mathrel{\preccurlyeq}}
\DeclareMathOperator*{\indlim}{\mathrm{ind\,lim}}
\DeclareMathOperator{\tr}    {\mathsf{tr}}
\newcommand{\tensor}[1][{}]           {\mathbin{\otimes_{\scriptscriptstyle{#1}}}}
\DeclarePairedDelimiter{\SP} {\langle}{\rangle}
\DeclareMathOperator{\spann} {\mathrm{span}}
\newcommand{\complete}[1]    {\widehat{#1}}
\DeclarePairedDelimiter{\norm}{\lVert}{\rVert}
\newcommand{\halbnorm}[1]       {\mathrm{#1}}
\newcommand{\algebra}[1]      {\mathcal{#1}}
\newcommand{\acts}            {\mathbin{\triangleright}}
\newcommand{\racts}           {\mathbin{\triangleleft}}
\newcommand{\basis}[1]{\mathsf{#1}}

\newcommand{\cpitensor}      {\mathbin{\hat{\otimes}_\pi}}
\newcommand{\ellOne}         {\ell^1}
\newcommand{\ellInfty}       {\ell^\infty}
\newcommand{\cNull}          {\mathfrak{c}^0}

\title{Nuclear Group Algebras for Finitely Generated Groups}

\author{\textbf{Michel Cahen}\thanks{\texttt{mcahen@ulb.ac.be}} \\[0.2cm]
  \begin{minipage}{8cm}
      \centering
      \small
      D\'epartement de Math\'ematique \\
      Universit\'e Libre de Bruxelles \\
      Campus Plaine, C. P. 218 \\
      Boulevard du Triomphe \\
      B-1050 Bruxelles \\
      Belgium
  \end{minipage}
  \vspace{0.4cm} \\
  \addtocounter{footnote}{2}
  \textbf{Simone Gutt}\thanks{\texttt{sgutt@ulb.ac.be}} \\[0.2cm]
  \begin{minipage}{6cm}
      \centering
      \small
      D\'epartement de Math\'ematique \\
      Universit\'e Libre de Bruxelles \\
      Campus Plaine, C. P. 218 \\
      Boulevard du Triomphe \\
      B-1050 Bruxelles \\
      Belgium
  \end{minipage}
  \begin{minipage}{6cm}
      \centering
      \small
      Universit\'e de Lorraine \`a Metz\\
      D\'epartement de Math\'ematique\\
      Ile du Saulcy \\
      F-57045~Metz Cedex 01\\
      France
  \end{minipage}
  \vspace{0.2cm} \\
  and
  \\[0.2cm]
  \textbf{Stefan Waldmann}\thanks{\texttt{Stefan.Waldmann@mathematik.uni-wuerzburg.de}}
  \\[0.2cm]
  \begin{minipage}{8cm}
      \centering
      \small
      Institut f\"ur Mathematik \\
      Universit\"at W\"urzburg \\
      Campus Hubland Nord \\
      Emil-Fischer-Stra\ss e 31 \\
      97074 W\"urzburg \\
      Germany
  \end{minipage}
}

\date{October 2016}

\begin{document}

%
%

\maketitle

%
%

\begin{abstract}
    We study completions of the group algebra of a finitely generated
    group and relate nuclearity of such a completion to growth
    properties of the group. This extends previous work of Jolissaint
    on nuclearity of rapidly decreasing functions on a finitely
    generated group to more general weights than polynomial
    decrease. The new group algebras and their duals are studied in
    detail and compared to other approaches. As application we discuss
    the convergence of the complete growth function introduced by
    Grigorchuk and Nagnibeda.
\end{abstract}

\newpage

%
%

\tableofcontents
\newpage

%
%

\section{Introduction}
\label{sec:Introduction}

In \cite{jolissaint:1990a, jolissaint:1989a} Jolissaint introduced and
studied a completion of the group algebra of a finitely generated
group resembling the construction of the Schwartz space of rapidly
decreasing functions: the coefficients decay faster than every
polynomial in the length of the group element, where one has to choose
a length functional for the group, e.g. based on the choice of
finitely many generators. This algebra turned out to be an extremely
useful tool to investigate other completions of the group algebra like
the (reduced) $C^*$-group algebra. Among many other results Jolissaint
showed that the Fr\'echet space of such rapidly decreasing functions is
nuclear iff the group has polynomial growth.

In this paper we generalize Jolissaint's constructions and results by
taking into account more general types of decays for the coefficients
in the group algebra.

To this end we consider growth functions $\sigma$ on the integers
which are either sub-multiplicative or almost sub-multiplicative (in
the sense of \eqref{eq:AlmostSubmultiplicative}). In both cases we
define, for a given choice of a length functional $L$ on the group and
a parameter $R \ge 0$, a weighted $\ell^1$-norm
\begin{equation*}
    \label{eq:ThisIsTheNorm}
    \norm{a}_{L, \sigma, R}
    =
    \sum_{g \in G} |a_g| \sigma(L(g))^R
\end{equation*}
for $a \in \mathbb{C}[G]$. The system of all such norms for $R \in [0,
\infty)$ yield an algebra topology on the group algebra
$\mathbb{C}[G]$ which turns out to be independent of $L$. The
completion is shown to be a Fr\'echet-Hopf $^*$-algebra
$\algebra{A}_\sigma(G)$. We characterize the nuclearity of the
completion $\mathcal{A}_\sigma(G)$ by comparing the growth of $\sigma$
to the growth of the group: we have nuclearity iff $\sigma$ grows at
least as fast as the group (in the sense of \eqref{eq:earlierDef}). The
basic idea of the proof is to realize the algebra
$\algebra{A}_\sigma(G)$ as a certain Köthe space such that one can
apply the easy Grothendieck-Pietsch criterion for nuclearity.  We also
characterize the dual and its rich bimodule structure quite explicitly
and give as a first application convergence results for the complete
growth function, introduced by Grigorchuk and Nagnibeda in
\cite{grigorchuk.nagnibeda:1997a}, as an element of the dual: again,
we get convergence and thus a holomorphic complete growth function iff
the growth of $\sigma$ is at least as fast as the one of the group.

The paper is organized as follows: in
Section~\ref{sec:NotationPreliminaries} we explain our notation and
recall some basic facts on finitely generated groups to set-up the
stage. In section \ref{sec:GrowthFunctions} we specify the kind of
growth functions which we are going to use and discuss many
examples. The general construction of the norms and completions of the
group algebra is presented in
Section~\ref{sec:GeneralConstruction}. In Section~\ref{sec:Comparison}
we compare our construction to yet another version from
\cite{beiser.waldmann:2014a}; understanding the construction presented
there was one of the motivations to start the present work.
Section~\ref{sec:Properties} contains the main result on the
nuclearity as well as a discussion of the canonical Schauder basis of
$\algebra{A}_\sigma(G)$. The dual and its inductive limit topology are
described in Section~\ref{sec:Dual}. Section~\ref{sec:CnullSpaces}
contains a detailed analysis of sub-bimodules of the dual which can be
seen as analogs of spaces of continuous functions vanishing (in a
weighted sense) at infinity. Finally, we prove the convergence
statements for the complete growth function in
Section~\ref{sec:TotalGrowth}. In a small appendix we collect some
background material in locally convex analysis of sequence spaces.

\noindent
\textbf{Acknowledgements:} This research was partially supported by
the IAP programme ``Dygest'' of the Belgium Science Policy Office.

%
%

\section{Notation and Preliminaries}
\label{sec:NotationPreliminaries}

In this section we will recall some basic notions from the theory of
discrete groups and their group algebras. Let $G$ be a discrete group,
i.e. we consider the group as a discrete topological space and hence
as a topological group. Later on, we will be interested in finitely
generated groups, but for the time being, this is not yet
important. We consider the \emph{group algebra} $\mathbb{C}[G]$
defined to be the complex vector space spanned by the set $G$ and
endowed with the multiplication inherited from $G$. More precisely, we
have the basis vectors $\basis{e}_g \in \mathbb{C}[G]$, where $g \in
G$, and the product of $a, b \in \mathbb{C}[G]$ is given by
\begin{equation}
    \label{eq:GroupAlgebraProduct}
    ab
    =
    \Bigg(\sum_{g \in G} a_g \basis{e}_g\Bigg)
    \Bigg(\sum_{h \in G} b_h \basis{e}_h\Bigg)
    =
    \sum_{g, h \in G} a_gb_h \basis{e}_{gh}
    =
    \sum_{h \in G}
    \Bigg(\sum_{g \in G} a_g b_{g^{-1}h}\Bigg)
    \basis{e}_h,
\end{equation}
where only finitely many of the coefficients $a_g$ and $b_h$ are
different from zero. We denote the coefficient functionals by
$\delta_g\colon \mathbb{C}[G] \ni a \mapsto a_g \in \mathbb{C}$.

The product \eqref{eq:GroupAlgebraProduct} makes the group algebra
$\mathbb{C}[G]$ a unital associative algebra where the unit element is
given by the group unit, i.e. by $\basis{e}_e$. Moreover,
$\mathbb{C}[G]$ becomes a \emph{cocommutative Hopf $^*$-algebra} by
the following structure maps: the $^*$-involution is given by $a^* =
\sum_{g \in G} \cc{a_g} \basis{e}_{g^{-1}}$, the antipode is $S(a) =
\sum_{g \in G} a_g \basis{e}_{g^{-1}}$, the counit is $\epsilon(a) =
\sum_{g \in G} a_g$, and the coproduct is $\Delta(a) = \sum_{g \in G}
a_g \basis{e}_g \tensor \basis{e}_g$. The \emph{group-like elements}
of $\mathbb{C}[G]$ with respect to this coproduct are then precisely
the group elements, i.e. the basis vectors $\basis{e}_g$ for $g \in
G$.  Note that the tensor product $\mathbb{C}[G] \tensor
\mathbb{C}[G]$ is canonically isomorphic as algebra to $\mathbb{C}[G
\times G]$ where $G \times G$ is endowed with the product group
structure.

In the following it will be advantageous to consider also the
\emph{formal power series} in $G$ which we denote by
$\mathbb{C}[[G]]$: this is the vector space of formal series $a =
\sum_{g \in G} a_g \basis{e}_g$ where now $a_g$ is unrestricted. Then
$\mathbb{C}[[G]]$ is no longer an algebra but the multiplication
\eqref{eq:GroupAlgebraProduct} extends to a \emph{bimodule structure}
for $\mathbb{C}[[G]]$ over $\mathbb{C}[G]$.

The formal power series $\mathbb{C}[[G]]$ can also be identified with
the algebraic dual of $\mathbb{C}[G]$ as follows. A linear functional
$\Phi\colon \mathbb{C}[G] \longrightarrow \mathbb{C}$ is determined by
its values $\varphi_g = \Phi(\basis{e}_g)$ on the basis
$\{\basis{e}_g\}_{g \in G}$. This defines a formal series $\varphi =
\sum_{g \in G} \varphi_g \basis{e}_g$ such that for all $a \in
\mathbb{C}[G]$ one has
\begin{equation}
    \label{eq:PairingDual}
    \Phi(a)
    = \sum_{g \in G} a_g \Phi(\basis{e}_g)
    = \sum_{g \in G} a_g \varphi_g
    = \SP{\varphi, a},
\end{equation}
where we take the last equation as definition for the bilinear pairing
$\SP{\argument, \argument}\colon \mathbb{C}[[G]] \times \mathbb{C}[G]
\longrightarrow \mathbb{C}$. Conversely, every $\varphi \in
\mathbb{C}[[G]]$ gives a linear functional $\Phi$ on $\mathbb{C}[G]$
by this equation. Using the coefficient functionals, the linear
functional $\Phi$ corresponding to $\varphi$ is just $\Phi = \sum_{g \in
  G} \varphi_g \delta_g$. The infinite series becomes finite after
application to an element in $\mathbb{C}[G]$.

While this is true for every vector space with basis, we have now the
additional bimodule structures of $\mathbb{C}[[G]]$ and of the dual:
clearly the dual $\algebra{A}^*$ of an algebra $\algebra{A}$ is always
a bimodule over the algebra by
\begin{equation}
    \label{eq:BimoduleDual}
    (a \acts \Phi \racts b)(c)
    =
    \Phi(bca)
\end{equation}
for $a, b, c \in \algebra{A}$ and $\Phi \in \algebra{A}^*$. We can
relate this canonical bimodule structure to the one of
$\mathbb{C}[[G]]$ as follows. For later use we introduce already here
the \emph{trace functional}
\begin{equation}
    \label{eq:Trace}
    \tr\colon \mathbb{C}[[G]] \ni a
    \; \mapsto \;
    a_e \in \mathbb{C}.
\end{equation}
With some abuse of notation we have $\tr = \delta_e$, but note that
$\delta_e$ is a linear functional on $\mathbb{C}[G]$ while $\tr$ is
defined on $\mathbb{C}[[G]]$.  It is now easy to see that $\tr$ is
indeed a trace in the sense that for $a \in \mathbb{C}[G]$ and $b \in
\mathbb{C}[[G]]$ we have
\begin{equation}
    \label{eq:TraceIsTrace}
    \tr(ab) = \tr(ba),
\end{equation}
where we use the above bimodule structure of $\mathbb{C}[[G]]$. In
particular, $\tr$ restricts to a trace on the algebra $\mathbb{C}[G]$
given by $\tr\at{\mathbb{C}[G]} = \delta_e$. For the bimodule
structure we note that
\begin{equation}
    \label{eq:BimoduleIsomorphism}
    \mathbb{C}[[G]] \ni b
    \; \mapsto \;
    \SP{S(b), \argument} = \tr(b \argument) \in (\mathbb{C}[G])^*
\end{equation}
is now an isomorphism of bimodules, where we use the canonical
bimodule structures for both sides. Note that the use of the antipode
is crucial here. As a consequence, we can embed the group algebra
$\mathbb{C}[G]$ into its algebraic dual $(\mathbb{C}[G])^*$ as a
\emph{bimodule}.

In the following we want to study completions of $\mathbb{C}[G]$ with
respect to various topologies. The easiest case is well-known and
given by the $\ell^1$-like group algebra: we want to extend the counit
$\epsilon$ to an as large as possible subspace of $\mathbb{C}[[G]]$
such that the formula $\epsilon(a) = \sum_{g \in G} a_g$ still makes
sense. Since we do not require a specific order of summation, we have
either unconditional convergence or divergence: but in $\mathbb{C}$
unconditional convergence is the same as absolute
convergence. Moreover, we know that at most countably many
contributions can be different from zero. This motivates the following
definition: for $a \in \mathbb{C}[[G]]$ one defines
\begin{equation}
    \label{eq:ellEinsSeminorm}
    \norm{a}_{\ell^1(G)} = \sum_{g \in G} |a_g|
    \in [0, +\infty],
\end{equation}
and sets $\ell^1(G) = \{a \in \mathbb{C}[[G]] \; | \;
\norm{a}_{\ell^1(G)} < \infty\}$. Then it is easy to see that
$\norm{\argument}_{\ell^1(G)}$ is a norm on $\ell^1(G)$ turning it
into a Banach space with $\mathbb{C}[G]$ being a dense
subspace. Moreover, $|\epsilon(a)| \le \norm{a}_{\ell^1(G)}$ shows
that $\epsilon$ is a continuous linear functional on $\ell^1(G)$.  The
space $\ell^1(G)$ is optimal with respect to this property: it is the
largest subspace to which $\epsilon$ can be extended by means of the
convergent series.  The basis vectors $\basis{e}_g$ of $\mathbb{C}[G]$
constitute an \emph{absolute Schauder basis} of $\ell^1(G)$, i.e. the
coefficient functionals $\delta_g$ are continuous and $a = \sum_{g \in
  G} \delta_g(a) \basis{e}_g$ converges absolutely for all $a \in
\ell^1(G)$ such that $\sum_{g \in G} |\delta_g(a)|
\norm{\basis{e}_g}_{\ell^1(G)}$ can be estimated by the norm: it
simply coincides with $\norm{a}_{\ell^1(G)}$. Note that for $a \in
\ell^1(G)$ at most countably many coefficients $\delta_g(a) = a_g$ are
different from zero.

The remaining algebraic features of $\mathbb{C}[G]$ extend
continuously to $\ell^1(G)$: the product is continuous making
$\ell^1(G)$ a Banach algebra, the $^*$-involution as well as the
antipode are continuous with operator norm $1$, and the coproduct
extends to $\ell^1(G)$ taking now values in $\ell^1(G) \cpitensor
\ell^1(G) = \ell^1(G \times G)$: here we have to use the completed
projective tensor product of the Banach space $\ell^1(G)$ with itself
and we identified the $\ell^1$-norm on $\mathbb{C}[G \times G]$ with
the tensor products of the $\ell^1$-norms according to
Lemma~\ref{lemma:TensorProductEllEinsSeminorms}. Then we have the
equality $\norm{\Delta(a)}_{\ell^1(G \times G)} =
\norm{a}_{\ell^1(G)}$ for all $a \in \ell^1(G)$, which proves the
continuity of $\Delta$. This way we obtain a \emph{Banach-Hopf
  $^*$-algebra} $\ell^1(G)$. It is now fairly easy to see that also in
this completion the group-like elements are simply given by the basis
vectors $\{\basis{e}_g\}_{g \in G}$ as before. Thus also $\ell^1(G)$
allows to recover $G$ as a group.

%
%

\section{Growth Functions}
\label{sec:GrowthFunctions}

In this section we recall the basic properties of the word length in a
finitely generated group and the resulting surface and volume growth
properties, see e.g. \cite[Chap.~VI]{delaharpe:2000a}. Moreover, we
will give a class of growth functions to which we want to compare the
growth properties of a finitely generated group.

Let $G$ be a finitely generated group with a finite set $S$ of
generators. The choice of $S$ then gives a \emph{word length} $L\colon
G \longrightarrow \mathbb{N}_0$ by taking the minimal number $L(g)$ of
generators needed to write $g \in G$ as a product of generators. By
convention, $L(e) = 0$ for the group unit $e \in G$ and we assume that
$S = S^{-1}$. Then $L$ enjoys the properties
\begin{equation}
    \label{eq:LProperties}
    L(g) = L(g^{-1})
    \quad
    \textrm{and}
    \quad
    L(gh) \le L(g) + L(h)
\end{equation}
for all $g, h \in G$. These properties allow to define a word metric
$d\colon G \times G \longrightarrow \mathbb{N}_0$ by setting
\begin{equation}
    \label{eq:WordMetric}
    d(g, h) = L(g^{-1}h)
\end{equation}
for $g, h \in G$. If $L'$ is another word length functional
corresponding to a different set of generators $S'$ then we have
constants $c, c' \in \mathbb{N}$ with $L(g) \le c L'(g)$ and $L'(g)
\le c'L(g)$ for all $g \in G$. This implies that the two word metrics
give the same coarse geometry: we have estimates $d(g, h) \le c d'(g,
h)$ and $d'(g, h) \le c' d(g, h)$ for all $g, h \in G$.

Using the word length one defines the \emph{surface growth} and the
\emph{volume growth} of the group $G$ by
\begin{equation}
    \label{eq:SurfaceGrowth}
    \sigma_G(n) = \#\{g \in G \; | \; L(g) = n\}
    \quad
    \textrm{and}
    \quad
    \beta_G(n) = \#\{g \in G \; | \; L(g) \le n\},
\end{equation}
where we suppress the dependence on $L$ in the notation. We have
$\sigma_G(n) \le \beta_G(n)$ and $\beta_G(n)$ grows at most
exponentially in $n$. Since we assume $G$ to be infinite, we have
$\sigma_G(n) \ge 1$ for all $n$ and $\beta_G$ is strictly
increasing. In particular $\beta_G(n) > n$. Moreover, we have the
following crucial properties
\begin{equation}
    \label{eq:SubmultSigmaBeta}
    \sigma_G(n+m) \le \sigma_G(n)\sigma_G(m)
    \quad
    \textrm{and}
    \quad
    \beta_G(n+m) \le \beta_G(n)\beta_G(m)
\end{equation}
for all $n, m \in \mathbb{N}_0$.

We shall now generalize the properties of the functions $\sigma_G$ and
$\beta_G$ in order to find a class of functions to which we can
compare them. Here we will not use the standard way of comparing
growth functions, as this is usually done in geometric group theory,
see \cite[Chap.~VI]{delaharpe:2000a}. Instead, we will need a
slightly coarser notion.

There are two aspects of importance here: the surface growth as well
as the volume growth are submultiplicative. The volume growth is also
monotonic.  In general, we call a map $\sigma\colon \mathbb{N}_0
\longrightarrow [0, \infty)$ \emph{submultiplicative} if
\begin{equation}
    \label{eq:Submultiplicative}
    \sigma(n+m) \le \sigma(n)\sigma(m)
\end{equation}
for all $n, m \in \mathbb{N}_0$ and \emph{almost submultiplicative} if
for every $\epsilon > 0$ there is a constant $c > 0$ such that
\begin{equation}
    \label{eq:AlmostSubmultiplicative}
    \sigma(n+m) \le  c \sigma(n)^{1+\epsilon}\sigma(m)^{1+\epsilon}
\end{equation}
for all $n, m \ge \mathbb{N}_0$. By induction, we conclude that for an
almost submultiplicative map $\sigma$ we get for every $\epsilon > 0$
and every $m \in \mathbb{N}$ a constant $c_{\epsilon, m} > 0$ with
\begin{equation}
    \label{eq:AlmostSubnm}
    \sigma(nm) \le c_{\epsilon, m} \sigma(n)^{m+\epsilon}
\end{equation}
for all $n \ge 0$. The almost submultiplicative functions will allow
for a slightly greater flexibility later on. Clearly, a
submultiplicative function is also almost submultiplicative.  It will
be important to consider only functions $\sigma$ which are non-zero
everywhere. Moreover, we will need that $\sigma$ is monotonic: here we
use monotonic in the sense of $\sigma(n) \le \sigma(n+1)$ and
emphasize the \emph{strict} monotonic case $\sigma(n) < \sigma(n+1)$
whenever it is actually needed.
\begin{definition}[Growth function]
    \label{definition:GrowthFunction}%
    A map $\sigma\colon \mathbb{N}_0 \longrightarrow [1, \infty)$ is
    called a growth function if it is monotonically increasing and
    unbounded.
\end{definition}
Later on, we will add some additional conditions concerning the speed
how fast $\sigma$ diverges to $+\infty$ for large arguments.
\begin{example}[Growth functions]
    \label{example:Growths}%
    The following examples of growth functions will be of particular
    importance in this work:
    \begin{examplelist}
    \item The volume growth $\beta_G$ of a finitely generated infinite
        group is a growth function which increases even strictly
        monotonically. Moreover, this growth function is
        submultiplicative. While the surface growth $\sigma_G$ is
        submultiplicative, it will not be a growth function in
        general, since there is no need for $\sigma_G$ to be
        unbounded: e.g. for $G = (\mathbb{Z}, +)$ the surface growth
        with respect to the standard generators $S = \{1, -1\}$ is
        constant $\sigma_{\mathbb{Z}}(n) = 2$.
    \item If $\sigma$ and $\tau$ are growth functions, then $x\sigma +
        y \tau$ is a growth function for all $x, y \ge 1$. Moreover,
        $\sigma\tau$ and $\sigma^d$ for $d > 0$ are growth
        functions.
    \item The map $\sigma(n) = 1+n$ is a submultiplicative growth
        function, called the \emph{linear growth function}. Note that
        we use $1+n$ in order to have non-zero values everywhere. Also
        $\sigma(n) = a_d n^d + \cdots + a_1 n + 1$ is a growth
        function whenever $a_d, \ldots, a_1 \ge 0$ and $a_d > 0$,
        which we call a \emph{polynomial growth} of degree $d$.
    \item Let $0 < \theta \le 1$ then $\sigma(n) = \exp(n^\theta)$ is
        a submultiplicative growth function. This is called the
        \emph{exponential growth} for $\theta = 1$ and a
        \emph{sub-exponential growth} for $\theta < 1$.
    \item The factorial $\sigma(n) = n!$ is a growth function, called
        the \emph{factorial growth}. In fact, it is not
        submultiplicative but almost submultiplicative: this follows
        since a binomial coefficient can always be dominated by an
        arbitrarily small power of a factorial. In fact, for all
        $\epsilon > 0$ we have $(n+m)! \le c_\epsilon
        (n!)^{1+\epsilon} (m!)^{1+\epsilon}$ where $c_\epsilon$ is
        chosen by e.g. $c_\epsilon = \sup_{n, m}
        \frac{2^{n+m}}{n!^\epsilon m!^\epsilon} < \infty$.
    \item With the convention $x! = \Gamma(x+1)$ for $x \in
        \mathbb{R}$, we can obtain a growth function for every $0 <
        \theta \le 1$ by setting $\sigma(n) = (n^\theta)!$. Since the
        $\Gamma$-function is (strictly) monotonically increasing for
        arguments $x \ge 2$, we obtain a strictly increasing growth
        function. We call these growth functions \emph{sub-factorial}.
    \item If $\tau\colon \mathbb{N}_0 \longrightarrow [0, \infty)$ is
        subadditive, i.e. $\tau(n + m) \le \tau(n) + \tau(m)$ and if
        $\sigma\colon [0, \infty) \longrightarrow [0, \infty)$ is
        submultiplicative then $\sigma \circ \tau$ is
        submultiplicative again. If $\sigma$ was only almost
        submultiplicative, then $\sigma \circ \tau$ is almost
        submultiplicative again. It follows that $\sigma$ from the
        previous part is an almost submultiplicative growth function
        for all $0 < \theta \le 1$. This can of course be checked
        directly as well. In fact, these (sub-) factorial growth
        functions constitute one of our main motivations.
    \end{examplelist}
\end{example}

We need to establish a reasonable notion of equivalence which measures
the type of growth of $\sigma$. Thus let $\sigma, \sigma'\colon
\mathbb{N}_0 \longrightarrow [1, \infty)$ be two maps. We define now a
relation $\earlier$ by $\sigma \earlier \sigma'$ if there are
constants $c, k \ge 1$ with
\begin{equation}
    \label{eq:earlierDef}
    \sigma(n) \le c \sigma'(cn)^k
\end{equation}
for all $n \in \mathbb{N}_0$. Though we are mainly interested in
growth functions, we will use this relation also for non-monotonic
functions.

It is now easy to see that this gives a transitive and reflexive
relation.  We will call two functions $\sigma$ and $\sigma'$
\emph{equivalent} if we have $\sigma \earlier \sigma'$ and $\sigma'
\earlier \sigma$. In this case we write $\sigma \sim \sigma'$. Note
that our notion of equivalence and $\earlier$ is slightly coarser than
the usual notion employed in geometric group theory, see
\cite[Chap.~VI]{delaharpe:2000a}. In fact, we get back to the more
common notion if we allow only for $k = 1$ in
\eqref{eq:earlierDef}. However, it will be important to work with this
less restrictive notion to compare growth functions. In particular, we
always have
\begin{equation}
    \label{eq:sigmaEquivPowerSigma}
    \sigma \sim \sigma^d
\end{equation}
for all $d \ge 1$.

It is now easy to see that the relation $\earlier$ is compatible with
pointwise convex combinations, products, and powers, i.e. we have for
$\sigma \earlier \tau$ and $\sigma' \earlier \tau'$ and $a, b, d \ge
1$ also
\begin{equation}
    \label{eq:earlierOfCombinations}
    \sigma\sigma' \earlier \tau\tau',
    \quad
    a\sigma + b\sigma' \earlier a\tau + b \tau',
    \quad
    \textrm{and}
    \quad
    \sigma^d \earlier \tau^d.
\end{equation}

Quite different from the usual notion for comparing growth functions,
the relation $\earlier$ is \emph{in-sensitive} to different polynomial
growths:
\begin{lemma}
    \label{lemma:PolynomialGrowthTheSame}%
    Let $\sigma$ be a growth function which grows at least as
    $(1+n)^\theta$ for some $\theta > 0$, i.e. for which we have
    $(1+n)^\theta \earlier \sigma$. Then for all polynomial growth
    functions $p$ we have $p\sigma \sim \sigma$. In particular, all
    polynomial growth functions are equivalent.
\end{lemma}
\begin{proof}
    This is just a repeated application of
    \eqref{eq:earlierOfCombinations}.
\end{proof}

On the other hand, the equivalence relation $\sim$ is sensitive to
different types of super-polynomial growths. In particular, for all $0
< \theta < \theta'$ we have for every polynomial growth $p$
\begin{equation}
    \label{eq:InequivalentGrowths}
    1 \earlier p \earlier \exp(n^\theta) \earlier \exp(n^{\theta'}),
\end{equation}
\emph{without} having equivalence. Moreover, in between the different
types of sub-exponential growths we have the sub-factorial ones
\begin{equation}
    \label{eq:FactorialBetweenExponential}
    \exp(n^\theta) \earlier (n^\theta)! \earlier \exp(n^{\theta'})
\end{equation}
for $0 < \theta < \theta'$, again without having equivalence.  This
gives us a wealth of inequivalent growth functions, both for the
submultiplicative and the almost submultiplicative case.

Finally, we recall that a submultiplicative map grows at most
exponentially, i.e. we have $\sigma(n) \le c^n$ with $c =
\max\{\sigma(0), \sigma(1)\}$. An almost submultiplicative map can
grow faster than exponentially like e.g. $\sigma(n) = n!$. However,
also the growth of an almost submultiplicative map is bounded e.g. by
$\sigma(n) \le c^{2^n}$ for some suitable $c > 0$.

%
%

\section{The General Construction}
\label{sec:GeneralConstruction}

In this section we present the main construction of various other
group algebras for $G$ based on the following simple idea: the
$\ell^1$-completion is the largest possible completion to have the
counit well-defined by means of the series expression. Hence we are
looking now for finer topologies than the $\ell^1$-topology which will
still give us the structure of locally convex algebras and, hopefully,
Hopf algebras. The topologies will be determined by means of seminorms
which are analogs of the $\ell^1$-norm but involve now different
weights for each basis vector $\basis{e}_g$ instead of the uniform
counting measure.

The main idea is to find weights depending on the word length of the
group element in such a way that we get seminorms allowing for a
continuous multiplication. To guarantee a continuous product we use a
submultiplicative or an almost submultiplicative growth function
$\sigma$. We fix a set of generators with corresponding word length
$L$ on $G$. Moreover, we consider $R \ge 0$ and set
\begin{equation}
    \label{eq:NormLsigmaR}
    \norm{a}_{L, \sigma, R}
    =
    \sum_{g \in G} |a_g| \sigma(L(g))^R
\end{equation}
for $a \in \mathbb{C}[[G]]$, allowing for the value $+\infty$. We set
\begin{equation}
    \label{eq:ALsigmaR}
    \ellOne_{L, \sigma, R}(G)
    =
    \left\{
        a \in \mathbb{C}[[G]]
        \; \Big| \;
        \norm{a}_{L, \sigma, R} < \infty
    \right\}.
\end{equation}
Then clearly $\ellOne_{L, \sigma, R}(G)$ is a subspace of
$\mathbb{C}[[G]]$ containing $\mathbb{C}[G]$ and $\norm{\argument}_{L,
  \sigma, R}$ is a norm on it, turning $\ellOne_{L, \sigma, R}(G)$
into a Banach space with $\mathbb{C}[G]$ being a dense
subspace. Moreover,
\begin{equation}
    \label{eq:ellEinsKleiner}
    \norm{a}_{\ell^1(G)} \le \norm{a}_{L, \sigma, R}
\end{equation}
for all $a \in \mathbb{C}[[G]]$ and hence $\ellOne_{L, \sigma, R}(G)
\subseteq \ell^1(G)$ is a continuous inclusion. Note that the property
$\sigma \ge 1$ for a growth function becomes crucial here. For $R = 0$
we have equality of the norms and hence of the Banach
spaces. Similarly, for $R \ge R'$ we have
\begin{equation}
    \label{eq:RgreaterRprime}
    \norm{a}_{L, \sigma, R'} \le \norm{a}_{L, \sigma, R},
\end{equation}
leading to the continuous inclusion $\ellOne_{L, \sigma, R}(G)
\subseteq \ellOne_{L, \sigma, R'}(G)$. Finally, it is clear that
the $^*$-involution, the counit, and the antipode are continuous with
respect to the norms $\norm{\argument}_{L, \sigma, R}$ for all choices
of $\sigma$ and $R$ as above. In fact, $\norm{a^*}_{L, \sigma, R} =
\norm{a}_{L, \sigma, R} = \norm{S(a)}_{L, \sigma, R}$ since we have
$L(g^{-1}) = L(g)$. The continuity of the product is slightly more
tricky: here we need to use the (almost) submultiplicativity of
$\sigma$:
\begin{lemma}
    \label{lemma:ProductContinuous}%
    Let $\sigma$ be a growth function and let $R \ge 0$.
    \begin{lemmalist}
    \item If $\sigma$ is submultiplicative then $\norm{ab}_{R, \sigma,
          L} \le \norm{a}_{L, \sigma, R} \norm{b}_{L, \sigma, R}$ for
        all $a, b \in \mathbb{C}[[G]]$ and hence $\ellOne_{L,
          \sigma, R}(G)$ is a Banach $^*$-algebra.
    \item If $\sigma$ is almost submultiplicative then for every
        $\epsilon > 0$ there is a constant $c > 0$ such that
        $\norm{ab}_{R, \sigma, L} \le c \norm{a}_{L, \sigma,
          R+\epsilon} \norm{b}_{L, \sigma, R+\epsilon}$ for all $a, b
        \in \mathbb{C}[[G]]$.
    \end{lemmalist}
\end{lemma}
\begin{proof}
    The proof consist in a straightforward estimate using
    $\sigma(L(h)) \le \sigma(L(g) + L(g^{-1}h))$ and the (almost)
    submultiplicativity of $\sigma$.
\end{proof}

It follows that in the almost submultiplicative case we do not get a
normed algebra. However, there is a rather easy way out: we have to
pass to a suitable projective limit which is possible thanks to the
continuous inclusions $\ellOne_{L, \sigma, R}(G) \subseteq \ellOne_{L,
  \sigma, R'}(G)$ for $R \ge R'$. Hence, for a fixed $R > 0$, we take
the intersection
\begin{equation}
    \label{eq:ALsigmaRMinus}
    \ellOne_{L, \sigma, R^-}(G)
    =
    \bigcap_{R \ge \epsilon > 0} \ellOne_{L, \sigma, R-\epsilon}(G),
\end{equation}
and equip it with the projective limit topology: more explicitly, this
is determined by \emph{all} the norms $\norm{\argument}_{L, \sigma,
  R-\epsilon}$ for $0 < \epsilon \le R$. Since then for all $\delta >
0$ we find an $0 < \epsilon < \delta$ and a constant $c > 0$ with
\begin{equation}
    \label{eq:LsigmaRminusdelta}
    \norm{ab}_{L, \sigma, R - \delta}
    \le
    c
    \norm{a}_{L, \sigma, R - \epsilon}
    \norm{b}_{L, \sigma, R - \epsilon},
\end{equation}
we conclude that the product is continuous in the projective limit
topology. Moreover, the estimate from \eqref{eq:RgreaterRprime} shows
that countably many of the seminorms are sufficient to determine the
topology:
\begin{corollary}
    \label{corollary:ARminus}%
    For an almost submultiplicative growth function $\sigma$ and $R >
    0$ the space $\ellOne_{L, \sigma, R^-}(G)$ is a unital Fr\'echet
    $^*$-algebra with $\mathbb{C}[G]$ being a dense subalgebra. We
    have continuous inclusions $\ellOne_{L, \sigma, R^-}(G) \subseteq
    \ellOne_{L, \sigma, R'{}^-}(G)$ whenever $R \ge R'$. The counit
    $\epsilon$ and the antipode $S$ are continuous.
\end{corollary}

We come now to the question whether the coproduct is still continuous
to make these new algebras also Banach-Hopf or Fr\'echet-Hopf algebras.
Here things turn out to be slightly more complicated. As before we use
the tensor product of the norm $\norm{\argument}_{L, \sigma, R}$ with
itself to get a norm on $\mathbb{C}[G \times G]$. From
Lemma~\ref{lemma:TensorProductEllEinsSeminorms} we then get
\begin{equation}
    \label{eq:normDeltaa}
    \norm{\Delta(a)}_{L, \sigma, R}
    =
    \sum_{g, h \in G}
    |\Delta(a)_{g, h}| \sigma(L(g))^R \sigma(L(h))^R
    =
    \sum_{g \in G}
    |a_g| \sigma(L(g))^{2R}
    =
    \norm{a}_{L, \sigma, 2R}
\end{equation}
for all $a \in \mathbb{C}[G]$. Thus the coproduct extends to a
continuous linear map
\begin{equation}
    \label{eq:CoproductContinuous}
    \Delta\colon
    \ellOne_{L, \sigma, R}(G)
    \longrightarrow
    \ellOne_{L, \sigma, 2R}(G)
    \cpitensor
    \ellOne_{L, \sigma, 2R}(G),
\end{equation}
but in general it will not be continuous for the (strictly) finer
topology determined by the norm $\norm{\argument}_{L, \sigma, R}$ on
the target side. Hence $\ellOne_{L, \sigma, R}(G)$ will not be a
Banach-Hopf algebra unless $R = 0$.

Nevertheless, the estimate (equality) in
\eqref{eq:CoproductContinuous} gives a hint how to establish the
continuity: we have to take again a projective limit, now $R
\longrightarrow \infty$, in order to have \emph{all} seminorms
$\norm{\argument}_{L, \sigma, R}$ at our disposal.
\begin{definition}[The algebra $\algebra{A}_\sigma$]
    \label{definition:Asigma}%
    Let $\sigma$ be an almost submultiplicative growth function. Then
    we define
    \begin{equation}
        \label{eq:AsigmaDef}
        \algebra{A}_\sigma(G)
        =
        \left\{
            a \in \mathbb{C}[[G]]
            \; \Big| \;
            \norm{a}_{L, \sigma, R} < \infty
            \textrm{ for all }
            R \ge 0
        \right\}
        =
        \bigcap_{R \ge 0}
        \ellOne_{L, \sigma, R}(G),
    \end{equation}
    equipped with the projective locally convex topology of all the
    seminorms $\norm{\argument}_{L, \sigma, R}$.
\end{definition}

We will now investigate some first properties of
$\algebra{A}_\sigma(G)$. The most important is that even though the
individual spaces $\ellOne_{L, \sigma, R}(G)$ depend on the choice of
the set of generators and hence on $L$, this dependence disappears for
$\algebra{A}_\sigma$: suppose that $L'$ is another word length
obtained from a different set of generators. Then we have $\sigma(mn)
\le \sigma(n)^m$ for all $n \in \mathbb{N}_0$ in the case $\sigma$ is
submultiplicative. For an almost submultiplicative $\sigma$ we get for
all $\epsilon > 0$ a constant $c > 0$ depending on $m$ with
$\sigma(mn) \le c \sigma(n)^{m + \epsilon}$ for all $n \in
\mathbb{N}_0$. In the submultiplicative case this gives
\begin{equation}
    \label{eq:normLprime}
    \norm{a}_{L, \sigma, R}
    =
    \sum_{g \in G}
    |a_g| \sigma(L(g))^R
    \le
    \sum_{g \in G}
    |a_g| \sigma(cL'(g))^R
    \le
    \sum_{g \in G}
    |a_g| \sigma(L'(g))^{cR}
    =
    \norm{a}_{L', \sigma, cR},
\end{equation}
where we used the monotonicity of $\sigma$ and the constant $c$ from
$L(g) \le c L'(g)$. Using \eqref{eq:AlmostSubnm}, in the almost
submultiplicative case we get for all $\epsilon > 0$ another constant
$c' > 0$ with
\begin{equation}
    \label{eq:OtherLAlmostSub}
    \norm{a}_{L, \sigma, R}
    \le
    \norm{a}_{L', \sigma, c'R + \epsilon}.
\end{equation}
In both cases we can get the following statement for the projective
limit $R \longrightarrow \infty$:
\begin{proposition}
    \label{proposition:AlgebraSigma}%
    Let $\sigma$ be an almost submultiplicative growth function.
    \begin{propositionlist}
    \item \label{item:AsigmaIndependent} The projective limit
        $\algebra{A}_\sigma(G)$ of the Banach spaces $\ellOne_{L,
          \sigma, R}(G)$ is independent of the chosen word length $L$.
    \item \label{item:AsigmaFrechetHopf} The group algebra
        $\mathbb{C}[G]$ is dense in $\algebra{A}_\sigma(G)$ and all
        algebraic structures are continuous yielding a Fr\'echet-Hopf
        $^*$-algebra structure on $\algebra{A}_\sigma(G)$.
    \item \label{item:AsigmaLMC} If $\sigma$ is submultiplicative then
        $\algebra{A}_\sigma(G)$ is a locally multiplicatively convex
        algebra: the norms $\norm{\argument}_{L, \sigma, R}$ are
        submultiplicative.
    \end{propositionlist}
\end{proposition}
Note that in the almost submultiplicative case the norms
$\norm{\argument}_{L, \sigma, R}$ are \emph{not} submultiplicative.
In both cases, passing to the projective limit solved two problems at
once: we get a continuous coproduct and independence on the choice of
the word length.

Since we know that $\algebra{A}_\sigma(G)$ is actually independent of
the choice of generators we can easily discuss the relation with group
morphisms. Let $G'$ be another finitely generated group and let
$\phi\colon G \longrightarrow G'$ be a group morphism. We can then
extend $\phi$ to a linear map $\phi\colon \mathbb{C}[[G]]
\longrightarrow \mathbb{C}[[G']]$ by
\begin{equation}
    \label{eq:phiExtended}
    \phi\left(\sum_{g \in G} a_g \basis{e}_g\right)
    =
    \sum_{g \in G} a_g \basis{e}_{\phi(g)}.
\end{equation}
Restricted to $\mathbb{C}[G]$ this extension takes values in
$\mathbb{C}[G']$ and yields a Hopf $^*$-algebra morphism. If $S$ is a
set of generators of $G$ as above we get a word length $L$ and we get
a set $\phi(S) \subseteq G'$ which we can extend to a set of
generators $S'$ of $G'$. This gives a word length $L'$ for $G'$ which
satisfies
\begin{equation}
    \label{eq:LHsmallerL}
    L'(\phi(g)) \le L(g)
\end{equation}
for $g \in G$. From this the monotonicity of $\sigma$ implies
immediately the estimate
\begin{equation}
    \label{eq:phiContinuous}
    \norm{\phi(a)}_{L', \sigma, R} \le \norm{a}_{L, \sigma, R}
\end{equation}
for all $a \in \mathbb{C}[[G]]$. In particular, restricted to the
Banach spaces $\ellOne_{L, \sigma, R}(G)$ and $\ellOne_{L', \sigma,
  R}(G)$ the map $\phi$ is a continuous linear operator with operator
norm $1$ since on the group units we get
$\norm{\phi(\basis{e}_e)}_{L', \sigma, R} = 1 = \norm{\basis{e}_e}_{L,
  \sigma, R}$. From this we get the following functoriality of the
construction of the algebra $\algebra{A}_\sigma(G)$:
\begin{theorem}[The algebra $\algebra{A}_\sigma(G)$]
    \label{theorem:FunctorAsigma}%
    Let $\sigma$ be an almost submultiplicative growth function. Then
    assigning $\algebra{A}_\sigma(G)$ to a finitely generated group
    $G$ gives a covariant functor from the category of finitely
    generated groups to the category of Fr\'echet-Hopf $^*$-algebras,
    where on morphisms we use the continuous extension as above. If in
    addition $\sigma$ is submultiplicative, then
    $\algebra{A}_\sigma(G)$ is in addition locally multiplicatively
    convex.
\end{theorem}
\begin{proof}
    The functoriality is clear and the remaining parts have been shown
    already.
\end{proof}
\begin{remark}
    \label{remark:LMC}%
    The submultiplicative case is of course nicer as here we have the
    \emph{entire holomorphic calculus} of locally multiplicatively
    convex algebras at hand, see e.g. \cite{michael:1952a}. In the
    almost submultiplicative case this is not true in general. In both
    situations we have many examples thanks to
    Example~\ref{example:Growths}.
\end{remark}
Since $\algebra{A}_\sigma(G) \subseteq \ell^1(G)$ is a Hopf subalgebra
containing $\mathbb{C}[G]$, the group-like elements of
$\algebra{A}_\sigma(G)$ are again just the basis vectors
$\{\basis{e}_g\}_{g \in G}$.

%
%

\section{Comparison with Previous Results}
\label{sec:Comparison}

In this section we relate our general construction of a Fr\'echet-Hopf
$^*$-algebra completion of the algebraic group algebra $\mathbb{C}[G]$
to some previous and well-known constructions.

In \cite{jolissaint:1989a, jolissaint:1990a} the polynomial growth
version was studied. In our notation this corresponds to the algebra
$\algebra{A}_\sigma(G)$ with $\sigma$ being the linear growth
$\sigma(n) = 1 + n$ or any other non-constant polynomial according to
Lemma~\ref{lemma:PolynomialGrowthTheSame}: the topology of
$\algebra{A}_\sigma$ will not depend on the particular choice of
$\sigma$ as long as it is polynomial as we shall see
Proposition~\ref{proposition:DifferentSigmas} in some larger
generality. Jolissaint then studied many properties of
$\algebra{A}_\sigma(G)$ and its completion including many applications
to $C^*$-algebraic versions of the group algebra.

In \cite{beiser.waldmann:2014a} a general construction of seminorms
for an algebra with countable vector space basis was proposed and
exemplified. Among many other examples, the group algebra of a
finitely generated group was discussed. We want to show that this
construction fits into our general approach for a specific choice of
$\sigma$. In fact, understanding the quite mysterious recursion scheme
of \cite{beiser.waldmann:2014a} better was one of the motivations for
the present paper.

Let us briefly recall the construction from
\cite{beiser.waldmann:2014a}: first we choose a function
$\mathsf{c}\colon G \longrightarrow \mathbb{R}^+$ with $\mathsf{c}(g)
= \mathsf{c}(g^{-1})$. Then out of this function one defines
recursively the numbers $h_{0, 0, g}(a) = |a_g|\mathsf{c}(g)$ and
\begin{equation}
    \label{eq:FunnyhNumbers}
    h_{m+1, 2\ell, k}(a)
    =
    \sum_{g \in G}
    h_{m, \ell, g}(a)^2
    \frac{\mathsf{c}(k)}{\mathsf{c}(g)\mathsf{c}(g^{-1}k)}
    \quad
    \textrm{and}
    \quad
    h_{m+1, 2\ell+1, k}(a)
    =
    h_{m+1, 2\ell, k^{-1}}(a)
\end{equation}
for $m \in \mathbb{N}$, $\ell = 0, \ldots, 2^m -1$ and $k \in G$. Then
the quantities
\begin{equation}
    \label{eq:TheFunnySeminorms}
    \norm{a}_{m, \ell, g}(a) = \sqrt[2^m]{h_{m, \ell, g}(a)}
\end{equation}
turn out to be seminorms which make the product of $\mathbb{C}[G]$
continuous \emph{provided} one can guarantee that
$\norm{a}_{m, \ell, g} < \infty$ is finite at all: this will depend on
the choice of the scaling function $\mathsf{c}$. While there are also
other possibilities, one successful choice in
\cite{beiser.waldmann:2014a} was to consider
$\mathsf{c}(g) = (L(g)!)^\rho$ for a fixed value $\rho > 0$. Then the
countable collection \eqref{eq:TheFunnySeminorms} of seminorms make
the group algebra a locally convex algebra which has a Fr\'echet algebra
completion denoted by $\algebra{A}_{L, \rho}(G)$.  To relate the two
topology we first show the following estimate:
\begin{lemma}
    \label{lemma:ComparisonEstimates}%
    Let $\sigma(n) = n!$ be the factorial growth and let $R, \rho > 0$
    be fixed. For $a \in \mathbb{C}[G]$ one has:
    \begin{lemmalist}
    \item \label{item:TechLemmaEstimateFunctionalag}
        $|a_g| \le (L(g)!)^{\rho(2^{-m}-1)} \norm{a}_{m+1, 0, e}$ for
        $m \in \mathbb{N}_0$.
    \item \label{item:OneDirectionEstimate} For large enough $m$ there
        is a constant $c > 0$ such that
        $\norm{a}_{L, \sigma, R} \le c \norm{a}_{m, 0, e}$ provided
        $R < \rho$.
    \item \label{item:OtherDirectionEstimate} For all
        $m \in \mathbb{N}_0$ and for all $0 < \epsilon < 4^{-m}\rho$
        there exists a constant $c_{m, \epsilon}$ such that for all
        $k \in G$ one has
        $h_{m, \ell, k}(a) \le c_{m, \epsilon} (L(k)!)^{4^m\epsilon}
        \norm{a}_{L, \sigma, \rho - \epsilon}^{2^m}$.
    \end{lemmalist}
\end{lemma}
\begin{proof}
    The first claim is \cite[Lem.~3.11]{beiser.waldmann:2014a} and
    follows directly from the definition and a little induction on
    $m$. For the second, we estimate
    \[
    \norm{a}_{L, \sigma, R}
    =
    \sum_{g \in G}
    |a_g| (L(g)!)^R
    \le
    \sum_{g \in G}
    (L(g)!)^R (L(g)!)^{{\rho(2^{-m}-1)}} \norm{a}_{m+1, 0, e}
    \]
    by the first part. Taking $m$ large enough, $R - \rho (1- 2^{-m})$
    becomes negative. But a negative power of the function
    $g \mapsto L(g)!$ is always summable over $G$, yielding the
    constant $c$ we need. For the third part we clearly find
    $c_{0, \epsilon}$ with the desired property. Thus we prove the
    estimate by induction on $m$. We have for
    $0 < \epsilon < 4^{-m-1} \rho < 4^{-m}\rho$ and hence
    \begin{align*}
        h_{m+1, 2\ell, k}(a)
        &\le
        \sum_{g \in G}
        c_{m, \epsilon}^2
        (L(g)!)^{4^m 2\epsilon}
        \norm{a}_{L, \sigma, \rho - \epsilon}^{2^{m+1}}
        \frac{(L(k)!)^\rho}{(L(g)!)^\rho (L(g^{-1}k)!)^\rho} \\
        &=
        c_{\epsilon, m}^2
        \norm{a}_{L, \sigma, \rho - \epsilon}^{2^{m+1}}
        \left(\sum_{(1)} + \sum_{(2)} + \sum_{(3)}\right)
        (L(g)!)^{4^m 2\epsilon}
        \frac{(L(k)!)^\rho}{(L(g)!)^\rho (L(g^{-1}k)!)^\rho},
    \end{align*}
    where we distinguish the following three cases for $g \in G$ in
    the summation: case ($1$) is $L(k) < L(g)$, case ($2$) is
    $L(k) \ge L(g)$ and $L(k) < L(g^{-1}k)$, and case ($3$) is
    $L(k) \ge L(g)$ and $L(k) \ge L(g^{-1}k)$. The first
    contribution can be estimated as
    \begin{align*}
        \sum_{(1)}
        (L(g)!)^{4^m 2\epsilon}
        \frac{(L(k)!)^\rho}{(L(g)!)^\rho (L(g^{-1}k)!)^\rho}
        &\le
        \sum_{(1)}
        (L(g)!)^{4^m 2\epsilon}
        \frac{1}{(L(g^{-1}k)!)^\rho} \\
        &\le
        \sum_{h \in G}
        \binom{L(k) + L(h)}{L(h)}^{4^m 2\epsilon}
        \frac{(L(k)!)^{4^m 2\epsilon}}{(L(h)!)^{\rho - 4^m 2\epsilon}}
        \\
        &\le
        {(L(k)!)^{4^m 2\epsilon}}
        2^{4^m 2 \epsilon L(k)}
        \sum_{h \in G}
        \frac{2^{4^m 2 \epsilon L(h)}}{(L(h)!)^{\rho - 4^m 2\epsilon}},
    \end{align*}
    using first the substitution $h = g^{-1}k$ and then the standard
    estimate for a binomial coefficient.  Now the last sum converges
    as we have a positive power of the factorial of the length in the
    denominator by assumption. This gives a constant independent of
    $k$ while the first terms can be estimated by another constant
    times $(L(k)!)^{4^{m+1}\epsilon}$. For the second we have
    \[
    \sum_{(2)}
    (L(g)!)^{4^m 2\epsilon}
    \frac{(L(k)!)^\rho}{(L(g)!)^\rho (L(g^{-1}k)!)^\rho}
    \le
    \sum_{g \in G} \frac{1}{(L(g)!)^{\rho - 4^m 2\epsilon}}
    <
    \infty,
    \]
    independently of $k$ since still $\rho - 4^m 2\epsilon > 0$. Thus
    the second contribution can be estimated by a constant times
    $(L(k)!)^{4^{m+1}\epsilon}$ directly. For the third contribution
    we get
    \begin{align*}
        \sum_{(3)}
        (L(g)!)^{4^m 2\epsilon}
        \frac{(L(k)!)^\rho}{(L(g)!)^\rho (L(g^{-1}k)!)^\rho}
        &\le
        \sum_{(3)}
        (L(g)!)^{4^m 2\epsilon}
        \binom{L(g) + L(g^{-1}k)}{L(g)}^\rho \\
        &\le
        2^{4^{m+1}\epsilon L(k)}
        \sum_{(3)}
        (L(k)!)^{4^m 2\epsilon},
    \end{align*}
    again by the usual estimate of the binomial coefficient. Now the
    condition $L(g) \le L(k)$ in ($3$) gives at most $(2N)^{L(k)}$
    terms and hence the whole sum grows at most exponentially in
    $L(k)$. Thus we can also estimate this last contribution by a
    constant times $(L(k)!)^{4^{m+1}\epsilon}$. Putting everything
    together gives the third claim.
\end{proof}
\begin{proposition}
    \label{proposition:CompareToCrypticSeminorms}%
    Let $\rho > 0$ and let $\sigma(n) = n!$ be the factorial
    growth. For the algebra $\algebra{A}_{L, \rho}(G)$ from
    \cite{beiser.waldmann:2014a} we have
    \begin{equation}
        \label{eq:FunnyAlgebraIsProjLimEllEins}
        \algebra{A}_{L, \rho}(G)
        =
        \ell^1_{L, \sigma, \rho^-}(G).
    \end{equation}
\end{proposition}
\begin{proof}
    Indeed, the mutual estimates from
    Lemma~\ref{lemma:ComparisonEstimates},
    \refitem{item:OneDirectionEstimate} and
    \refitem{item:OtherDirectionEstimate}, directly show that the two
    locally convex topologies on $\mathbb{C}[G]$ coincide. Hence the
    corresponding completions coincide as well.
\end{proof}

This finally puts the rather ad-hoc and mysterious construction of
\cite{beiser.waldmann:2014a} into a larger context. It also emphasizes
the importance of the factorial growth. Note that the algebra
$\algebra{A}_{L, \rho}(G)$ will in general depend on the choice of the
length functional $L$: only after the projective limit $\rho
\longrightarrow \infty$ we get independence. Note also that
\cite{beiser.waldmann:2014a} was a much more general construction
leading to non-trivial results far beyond group algebras. In fact, the
construction of a locally convex topology for star product algebras
was the primary goal of that work.

%
%

\section{Properties of $\algebra{A}_\sigma$}
\label{sec:Properties}

In this section we shall investigate general properties of the functor
$\algebra{A}_\sigma$ and the resulting algebras: we are interested in
the role of $\sigma$, the nuclearity properties, and in the
topological bases of the resulting algebras.

We begin with the dependence on $\sigma$.  Here we first observe the
following statement:
\begin{lemma}
    \label{lemma:GrowthGivesInequalities}%
    Let $\sigma \earlier \sigma'$ with constants $c, k \ge 1$ as in
    \eqref{eq:earlierDef}. Without restriction we assume $c, k \in
    \mathbb{N}$.
    \begin{lemmalist}
    \item If $\sigma$ and $\sigma'$ are submultiplicative then for all
        $R > 0$ we have $\norm{a}_{L, \sigma, R} \le c^R \norm{a}_{L,
          \sigma', ckR}$.
    \item If $\sigma$ and $\sigma'$ are almost submultiplicative then
        we have for all $R > 0$ and $\epsilon > 0$ a constant $c'$
        with $\norm{a}_{L, \sigma, R} \le c' \norm{a}_{L, \sigma', ckR
          + \epsilon}$.
    \end{lemmalist}
\end{lemma}
\begin{proof}
    Let $\sigma$ be almost submultiplicative and let $c_{\epsilon, c}
    > 0$ be such that $\sigma(nm) \le c_{\epsilon, c}
    \sigma(n)^{m+\epsilon}$, see \eqref{eq:AlmostSubnm}. Then
    \[
    \norm{a}_{L, \sigma, R}
    =
    \sum_{g \in G} |a_g| \sigma(L(g))^R
    \le
    c^R \sum_{g \in G} |a_g| \sigma'(cL(g))^{kR}
    \le
    c^R c_{\epsilon, c}^{kR} \sum_{g \in G}
    |a_g| \sigma'(L(g))^{ckR + \epsilon kR}.
    \]
    Rescaling $\epsilon$ and renaming the constants gives the result
    for the almost submultiplicative case, the submultiplicative case
    is even easier.
\end{proof}

This estimate for the norms gives immediately continuous inclusions
\begin{equation}
    \label{eq:AsigmaIntoAOtherSigma}
    \ellOne_{L, \sigma', ckR}(G)
    \longrightarrow
    \ellOne_{L, \sigma, R}(G)
\end{equation}
in the submultiplicative case and analogously
\begin{equation}
    \label{eq:AlmostAsigmaIntoAOtherSigma}
    \ellOne_{L, \sigma', ckR+\epsilon}(G)
    \longrightarrow
    \ellOne_{L, \sigma, R}(G)
\end{equation}
in the almost submultiplicative case. Since we are interested in the
projective limit $R \longrightarrow \infty$, the differences between
the submultiplicative and the almost submultiplicative case disappear
and we get the following statement:
\begin{proposition}
    \label{proposition:DifferentSigmas}%
    Let $\sigma$ and $\sigma'$ be almost submultiplicative growth
    functions with $\sigma \earlier \sigma'$. Then we have the
    continuous inclusion $\algebra{A}_{\sigma'}(G) \longrightarrow
    \algebra{A}_\sigma(G)$. In particular, for $\sigma \sim \sigma'$
    we have
    \begin{equation}
        \label{eq:SameAlgebraDifferentSigma}
        \algebra{A}_\sigma(G) = \algebra{A}_{\sigma'}(G).
    \end{equation}
\end{proposition}

In a next step we consider the basis $\{\basis{e}_g\}_{g \in G}$ of
the algebraic group algebra. This turns out to be an absolute Schauder
basis of $\algebra{A}_\sigma(G)$: First we recall that the topology of
$\algebra{A}_\sigma(G)$ is \emph{finer} than the $\ell^1$-topology
where we have constant weights in the definition of the $\ell^1$-norm
instead of growing ones. Hence the evaluation functionals
$\delta_g\colon \ell^1(G) \longrightarrow \mathbb{C}$ restrict to
continuous linear functionals on $\algebra{A}_\sigma(G)$.  Moreover,
we note that
\begin{equation}
    \label{eq:NormOnBasis}
    \norm{\basis{e}_g}_{L, \sigma, R} = \sigma(L(g))^R
\end{equation}
for all $g \in G$ by the very definition of the norms
$\norm{\argument}_{L, \sigma, R}$.
\begin{proposition}
    \label{proposition:SchauderBasis}%
    Let $\sigma$ be an almost submultiplicative growth function. Then
    the vectors $\{\basis{e}_g\}_{g \in G}$ together with the
    evaluation functionals $\{\delta_g\}_{g \in G}$ form an absolute
    Schauder basis of $\algebra{A}_\sigma(G)$.
\end{proposition}
\begin{proof}
    Let $a \in \algebra{A}_\sigma(G)$ and let $\norm{\argument}_{L,
      \sigma, R}$ be one of the defining seminorms. Then we have
    \[
    \sum_{g \in G} |\delta_g(a)| \norm{\basis{e}_g}_{L, \sigma, R}
    =
    \sum_{g \in G} |a_g| \sigma(L(g))^R
    =
    \norm{a}_{L, \sigma, R},
    \]
    which is the estimate (even equality) we need to show that we have
    an absolute Schauder basis, see \eqref{eq:WhatIsSchauderBasis},
    since clearly $a = \sum_{g \in G} a_g \basis{e}_g$ converges
    (unconditionally) in the topology of $\algebra{A}_\sigma(G)$.
\end{proof}
\begin{remark}
    \label{remark:KoetheSpace}%
    Having an absolute Schauder basis is a rather strong property of a
    locally convex space. It follows that it is a \emph{Köthe space},
    see e.g.  \cite[Sect.~1.7.F]{jarchow:1981a} or
    \cite[§27]{meise.vogt:1992a}. In our case, the \emph{Köthe matrix}
    of $\algebra{A}_\sigma(G)$ is given by the matrix
    \begin{equation}
        K(\algebra{A}_\sigma(G))
        =
        \left(
            \norm{\basis{e}_g}_{L, \sigma, n}
        \right)_{g \in G, n \in \mathbb{N}}
        =
        \left(
            \sigma(L(g))^n
        \right)_{g \in G, n \in \mathbb{N}},
    \end{equation}
    viewed as $G \times \mathbb{N}$ matrix. Many locally convex
    properties of $\algebra{A}_\sigma(G)$ are now encoded in this
    matrix.
\end{remark}

We come now to the main property of $\algebra{A}_\sigma(G)$:
nuclearity. Since we have a Fr\'echet topology, it is reasonable to ask
whether it is nuclear or not.
\begin{lemma}
    \label{lemma:VolumeSurfaceGrowthSmallerSigma}%
    For a growth function $\sigma$ with $(1+n) \earlier \sigma$ we
    have $\sigma_G \earlier \sigma$ iff $\beta_G \earlier \sigma$.
\end{lemma}
\begin{proof}
    Since $\sigma_G(n) \le \beta_G(n)$, one direction is trivial. Thus
    let $\sigma_G \earlier \sigma$ and hence $\sigma_G(n) \le c
    \sigma(cn)^k$ for some $c, k \ge 1$. Then
    \[
    \beta_G(n)
    =
    \sum_{\ell=0}^n \sigma_G(\ell)
    \le
    \sum_{\ell=0}^n c\sigma(c\ell)^k
    \le
    \sum_{\ell=0}^n c\sigma(cn)^k
    =
    c(1+n) \sigma(cn)^k,
    \]
    from which we get $\beta_G \earlier (1+n)\sigma$ But for $(1+n)
    \earlier \sigma$ we have $(1+n)\sigma \earlier \sigma^2 \sim
    \sigma$ by \eqref{eq:sigmaEquivPowerSigma} and
    \eqref{eq:earlierOfCombinations}.
\end{proof}

It is quite reasonable to consider growth functions $\sigma$ which
grow at least polynomially: a finitely generated infinite group grows
at least linearly. If $(1 + n) \earlier \sigma$ we get
\begin{equation}
    \label{eq:InvSigmaSummable}
    \frac{1}{\sigma} \in \ell^p
\end{equation}
for some $p \ge 1$ by using the very definition of the relation
$\earlier$. For such growth functions, the following theorem gives a
complete answer to the question of nuclearity:
\begin{theorem}[Nuclearity of $\algebra{A}_\sigma(G)$]
    \label{theorem:Nuclear}%
    Let $(1+n) \earlier \sigma$ be an almost submultiplicative growth
    function. Then $\algebra{A}_\sigma(G)$ is nuclear iff $\beta_G
    \earlier \sigma$.
\end{theorem}
\begin{proof}
    There are many ways to define nuclearity. Here we rely on the
    Grothendieck-Pietsch criterion for sequence spaces with absolute
    Schauder bases, see e.g. \cite[Thm.~28.15]{meise.vogt:1992a} or
    \cite[Sect.~6.1]{pietsch:1972a}: since the topology is defined by
    means of weighted $\ell^1$-norms the Grothendieck-Pietsch
    criterion says that $\algebra{A}_\sigma(G)$ is nuclear iff for
    every continuous seminorm $\norm{\argument}_{L, \sigma, R}$ we
    find another continuous seminorm $\norm{\argument}_{L, \sigma,
      R'}$ from our defining system of seminorms such that the ratio
    \[
    \frac{\norm{\basis{e}_g}_{L, \sigma, R}}
    {\norm{\basis{e}_g}_{L, \sigma, R'}}
    =
    \frac{1}{\sigma(L(g))^{R' - R}}
    \]
    is summable over $G$. This is clearly equivalent to say that we
    find a $\rho > 0$ such that $\frac{1}{\sigma(L(g))^\rho}$ is
    summable over $G$, by taking $R'$ large enough.  Now first we note
    that $\beta_G \earlier \sigma$ iff $\sigma_G \earlier \sigma$ by
    Lemma~\ref{lemma:VolumeSurfaceGrowthSmallerSigma}. Hence we can
    use the more adapted surface growth $\sigma_G$ instead. We first
    assume $\sigma_G \earlier \sigma$ and choose $c, k > 0$ with
    $\sigma_G(n) \le c \sigma(cn)^k$ for all $n \in
    \mathbb{N}_0$. Since $\sigma$ is almost submultiplicative, we find
    for every $\epsilon > 0$ a $c' > 0$ such that
    \[
    \sum_{g \in G}
    \frac{1}{\sigma(L(g))^\rho}
    =
    \sum_{n=0}^\infty
    \frac{\sigma_G(n)}{\sigma(n)^\rho}
    \le
    \sum_{n=0}^\infty
    \frac{c \sigma(cn)^k}{\sigma(n)^\rho}
    \le
    c' \sum_{n=0}^\infty
    \frac{\sigma(n)^{ck+\epsilon}}{\sigma(n)^\rho},
    \]
    according to \eqref{eq:AlmostSubnm}. If $\sigma$ is even
    submultiplicative, it suffices to take $\epsilon = 0$ and $c' =
    1$. Since by assumption $\frac{1}{\sigma} \in \ell^p$ for some $p
    \ge 1$, we can take $\rho$ large enough to see that the right hand
    side is summable. This shows that $\algebra{A}_\sigma(G)$ is
    nuclear in this case. Conversely, assume $\sigma_G \earlier
    \sigma$ does \emph{not} hold. By taking $c = 1$ we find for all $k
    \in \mathbb{N}$ an index $n_k$ such that $\sigma_G(n_k) >
    \sigma(n_k)^k$. Now fix $\rho > 0$ then for all $k \ge \rho$ we
    have a $n_k$ with $\sigma_G(n_k) > \sigma(n_k)^\rho$. Hence the
    quotient sequence $\frac{\sigma_G(n)}{\sigma(n)^\rho}$ can not be
    a zero sequence, let alone summable. Thus $\algebra{A}_\sigma(G)$
    is not nuclear.
\end{proof}
\begin{remark}
    \label{remark:NuclearityCondition}%
    The condition $\sigma_G \earlier \sigma$, which is equivalent to
    $\beta_G \earlier \sigma$ according to
    Lemma~\ref{lemma:VolumeSurfaceGrowthSmallerSigma}, can also be
    understood as follows: for an at least polynomially growing
    $\sigma$, we have $\sigma_G \earlier \sigma$ iff
    \begin{equation}
        \label{eq:SummableSigma}
        \frac{1}{L^*\sigma} \in \ell^p(G)
    \end{equation}
    for some $p \ge 1$, where $L$ is the word length and $L^*\sigma =
    \sigma \circ L$ is the pull-back.
\end{remark}
\begin{corollary}
    \label{corollary:AlwaysNuclear}%
    If $\sigma$ is an almost submultiplicative growth growing at least
    exponentially, then $\algebra{A}_\sigma(G)$ is always nuclear.
\end{corollary}
\begin{corollary}
    \label{corollary:PolynomialGrowthOfG}%
    If $G$ is a group with polynomial growth then
    $\algebra{A}_\sigma(G)$ is always nuclear.
\end{corollary}
The following particular case was treated by Jolissaint
\cite[Thm.~3.1.7]{jolissaint:1990a}:
\begin{corollary}
    \label{corollary:PolynomialGrowth}%
    Let $\sigma$ be the polynomial growth. Then
    $\algebra{A}_\sigma(G)$ is nuclear iff $G$ has polynomial growth.
\end{corollary}

Yet another aspect of the condition $\beta_G \earlier \sigma$ is the
following: in our definition of the norm $\norm{\argument}_{L, \sigma,
  R}$ we have chosen a $\ell^1$-version since this will yield good
algebraic properties. However, for many other reasons other
$\ell^p$-version might be useful as well: the $p = \infty$ case will
be needed for the dual in the next section, the $\ell^2$-version has
important applications when it comes to comparison with
$C^*$-algebraic versions of the group algebra as discussed e.g. in
\cite{jolissaint:1989a}. For $p \ge 1$ as well as for $p = \infty$ we
define
\begin{equation}
    \label{eq:WeightedSupNorm}
    \norm{a}_{\ell^p, L, \sigma, R}
    =
    \sqrt[p]{
      \sum_{g \in G}
      |a_g|^p \sigma(L(g))^R
    }
    \quad
    \textrm{and}
    \quad
    \norm{a}_{\ell^\infty, L, \sigma, R}
    =
    \sup_{g \in G}
    |a_g| \sigma(L(g))^R
    \in [0, \infty]
\end{equation}
for $a \in \mathbb{C}[[G]]$ and set
\begin{equation}
    \label{eq:ellInftyG}
    \ell^p_{L, \sigma, R}(G)
    =
    \left\{
        a \in \mathbb{C}[[G]]
        \; \Big| \;
        \norm{a}_{\ell^p, L, \sigma, R}
        < \infty
    \right\}.
\end{equation}
Clearly, we obtain Banach spaces and continuous inclusions
\begin{equation}
    \label{eq:InclusionsEllp}
    \ellOne_{L, \sigma, R}(G)
    \subseteq
    \ell^p_{L, \sigma, R}(G)
    \subseteq
    \ell^q_{L, \sigma, R}(G)
    \subseteq
    \ellInfty_{L, \sigma, R}(G)
\end{equation}
for $p \le q$ according to the obvious norm estimates
\begin{equation}
    \label{eq:NormEstimatesDifferentEllp}
    \norm{a}_{\ell^\infty, L, \sigma, R}
    \le
    \norm{a}_{\ell^q, L, \sigma, R}
    \le
    \norm{a}_{\ell^p, L, \sigma, R}
    \le
    \norm{a}_{\ell^1, L, \sigma, R}
    =
    \norm{a}_{L, \sigma, R}.
\end{equation}
In general, all these topologies will be quite different. However, the
nuclearity condition $\beta_G \earlier \sigma$ helps to simplify their
mutual relations drastically:
\begin{proposition}
    \label{proposition:ellpEquivalentNorms}%
    Let $\sigma$ be an almost submultiplicative growth function with
    $\beta_G \earlier \sigma$. Then for every $R \ge 0$ there exists a
    $R' > R$ and a constant $c$ such that
    \begin{equation}
        \label{eq:ellOneEstimateByEllInfty}
         \norm{a}_{L, \sigma, R}
         \le
         c \norm{a}_{\ell^\infty, L, \sigma, R'}
    \end{equation}
    for all $a \in \mathbb{C}[[G]]$.
\end{proposition}
\begin{proof}
    This is essentially the same argument as needed for
    Theorem~\ref{theorem:Nuclear}. Let $a \in \mathbb{C}[[G]]$. Then
    \[
    \norm{a}_{L, \sigma, R}
    =
    \sum_{g \in G} |a_g| \sigma(L(g))^{R'}
    \frac{1}{\sigma(L(g))^{R'-R}}
    \le
    \norm{a}_{\ell^\infty, L, \sigma, R'}
    \sum_{g \in G}
    \frac{1}{\sigma(L(g))^{R'-R}}
    \]
    as an inequality in $[0, +\infty]$. If $\beta_G \earlier \sigma$
    we know that for large enough $R'$ the remaining series becomes
    finite, hence giving the constant $c$ needed for
    \eqref{eq:ellOneEstimateByEllInfty}.
\end{proof}

Together with \eqref{eq:NormEstimatesDifferentEllp} it follows that in
the projective limit $R \longrightarrow \infty$ it is irrelevant which
$\ell^p$-version we use to define the topology on
$\algebra{A}_\sigma(G)$: all yield the same nuclear Fr\'echet topology.

%
%

\section{The Dual}
\label{sec:Dual}

In this section we shall investigate the dual of
$\algebra{A}_\sigma(G)$ and determine some of its basic
properties. Since we have an (absolute) Schauder basis, the dual can
be described by the standard techniques for Köthe spaces. The basic
idea is to determine the duals of the spaces $\ellOne_{L, \sigma,
  R}(G)$ and pass to the projective limit afterwards: under lucky
circumstances this will be an inductive limit for the duals.

If $\Phi\colon \algebra{A}_\sigma(G) \longrightarrow \mathbb{C}$ is a
continuous linear functional, then it is uniquely determined by its
values on the absolute Schauder basis $\{\basis{e}_g\}_{g \in
  G}$. Hence on $a \in \mathbb{C}[G] \subseteq \algebra{A}_\sigma(G)$
it is given by the pairing $\Phi(a) = \SP{\varphi, a}$ with a unique
$\varphi \in \mathbb{C}[[G]]$. However, not all $\varphi$ give a
continuous extension of the pairing to a continuous linear functional
on the whole algebra $\algebra{A}_\sigma(G)$. Here we need to estimate
the growth of the coefficients by means of the weighted supremum norm
\begin{equation}
    \label{eq:TheImportantSupNorm}
    \norm{b}_{\ell^\infty, L, \sigma, -R}
    =
    \sup_{g \in G}
    \frac{|b_g|}{\sigma(L(g))^R}.
\end{equation}
Note that we need to take a \emph{negative} value $-R \le 0$ of the
parameter in the definition of the supremum norm
\eqref{eq:WeightedSupNorm}.

From the general statements in Appendix~\ref{sec:LocallyConvex} we
know that $\ellInfty_{L, \sigma, -R}(G)$ is isomorphic to the
topological dual of the Banach space $\ellOne_{L, \sigma, R}(G)$.  We
shall now make this isomorphism explicit such that the bimodule
structure is continuous as well. To this end, we first note that the
antipode $S$ is a continuous operator on $\ellInfty_{L, \sigma,
  -R}(G)$. In fact, we have
\begin{equation}
    \label{eq:AntipodebIsometric}
    \norm{S(b)}_{\ell^\infty, L, \sigma, -R}
    =
    \norm{b}_{\ell^\infty, L, \sigma, -R}
\end{equation}
since $L(g^{-1}) = L(g)$.

The continuity of the bimodule structure is controlled by the
following estimates, which we have to formulate for the
submultiplicative and the almost submultiplicative case separately:
\begin{lemma}
    \label{lemma:BimodContinuity}%
    Let $a \in \mathbb{C}[G]$ and $b \in \mathbb{C}[[G]]$ and let $h
    \in G$. Moreover, let $R \ge 0$.
    \begin{lemmalist}
    \item \label{item:BimoduleSubCase} For a submultiplicative growth
        function $\sigma$ we have
        \begin{equation}
            \label{eq:SubBimoduleCase}
            |(ab)_h|
            \le
            \sigma(L(h))^R
            \norm{a}_{L, \sigma, R}
            \norm{b}_{\ell^\infty, L, \sigma, -R}.
        \end{equation}
        The same estimate holds for $(ba)_h$.
    \item \label{item:BimoduleSubCase} For an almost submultiplicative
        growth function $\sigma$ we have for all $\epsilon > 0$ a
        constant $c > 0$ such that
        \begin{equation}
            \label{eq:AlmostSubBimoduleCase}
            |(ab)_h|
            \le
            c
            \sigma(L(h))^{R+\epsilon}
            \norm{a}_{L, \sigma, R + \epsilon}
            \norm{b}_{\ell^\infty, L, \sigma, -R}.
        \end{equation}
        The same estimate holds for $(ba)_h$.
    \end{lemmalist}
\end{lemma}
\begin{proof}
    We prove the second statement, the first is analogous but
    simpler. We fix $\epsilon > 0$ and choose $c \ge 1$ such that
    $\sigma(n + m) \le c
    \sigma(n)^{1+\epsilon}\sigma(m)^{1+\epsilon}$. Then we get
    \begin{align*}
        |(ab)_h|
        &=
        \left|
            \sum_{g \in G} a_g \sigma(L(g^{-1}h))^R
            \frac{b_{g^{-1}h}}{\sigma(L(g^{-1}h))^R}
        \right| \\
        &\le
        \sum_{g \in G} |a_g|
        c
        \sigma(L(g))^{R (1+\epsilon)}
        \sigma(L(h))^{R (1+\epsilon)}
        \frac{|b_{g^{-1}h}|}{\sigma(L(g^{-1}h))^R} \\
        &\le
        c \sigma(L(h))^{R(1+\epsilon)}
        \norm{a}_{L, \sigma, R(1+\epsilon)}
        \sup_{k \in G}
        \frac{|b_k|}{\sigma(L(k))^R}.
    \end{align*}
    Since this can be done for all $\epsilon > 0$, a simple rescaling
    gives \eqref{eq:AlmostSubBimoduleCase}.
\end{proof}

Remarkably, for $h = e$ the estimate \eqref{eq:AlmostSubBimoduleCase}
can be improved to
\begin{equation}
    \label{eq:abAteAlmostSub}
    |(ab)_e| \le
    \norm{a}_{L, \sigma, R}
    \norm{b}_{\ell^\infty, L, \sigma, -R}
\end{equation}
also for the almost submultiplicative case by a direct estimate. For
general $h$ however, we get only the weaker estimate from
\eqref{eq:AlmostSubBimoduleCase}.  Based on these estimates the
following proposition is now a simple consequence of the algebraic
properties of the trace functional $\tr$ outlined in
Section~\ref{sec:NotationPreliminaries}.
\begin{proposition}
    \label{proposition:DualALsigmaR}%
    Let $R \ge 0$ and let $\sigma$ be an almost submultiplicative
    growth function. Then the map
    \begin{equation}
        \label{eq:DualIslInfty}
        \ellInfty_{L, \sigma, -R}(G) \ni b
        \; \mapsto \;
        \SP{S(b), \argument} = \tr(b \argument)
        \in \ellOne_{L, \sigma, R}(G)'
    \end{equation}
    is an isometric isomorphism of Banach spaces.
\end{proposition}
\begin{proof}
    From the general results on the dual of weighted $\ell^1$-spaces
    as in the Appendix~\ref{sec:LocallyConvex} and the algebraic
    considerations from Section~\ref{sec:NotationPreliminaries} this
    follows now from the estimate in \eqref{eq:abAteAlmostSub}.
\end{proof}

While \eqref{eq:DualIslInfty} is always an isomorphism of Banach
spaces, the situation for the bimodule structures is slightly more
complicated. Here we have to distinguish the submultiplicative and the
almost submultiplicative case. For a submultiplicative $\sigma$ we
know that $\ellOne_{L, \sigma, R}(G)$ is a Banach algebra and
hence its dual is a Banach bimodule over it in a canonical way. For
the almost submultiplicative case we first have to pass from
$\ellOne_{L, \sigma, R}(G)$ to the projective limit
$\ellOne_{L, \sigma, R^-}(G)$ in order to get an algebra at
all. Ultimately, we want to consider $R \longrightarrow \infty$ for
both cases.

To this end we compare the different $\ell^\infty$-norms: obviously we
have
\begin{equation}
    \label{eq:InftyNormsDifferentR}
    \norm{b}_{\ell^\infty, L, \sigma, -R}
    \le
    \norm{b}_{\ell^\infty, L, \sigma, -R'}
   \end{equation}
for $R \ge R' \ge 0$. Moreover, for $R = 0$ we have the usual supremum
norm. This results in continuous inclusions
\begin{equation}
    \label{eq:InclusionsAintoEllInfty}
    \ellOne_{L, \sigma, R}
    \longrightarrow
    \ell^1(G)
    \longrightarrow
    \ellInfty(G)
    \longrightarrow
    \ellInfty_{L, \sigma, -R'}(G)
    \longrightarrow
    \ellInfty_{L, \sigma, -R}(G)
\end{equation}
for $R \ge R' \ge 0$.
\begin{proposition}
    \label{proposition:BimoduleSubCase}%
    Let $\sigma$ be a submultiplicative growth function and let $R \ge
    0$.
    \begin{propositionlist}
    \item For $a \in \mathbb{C}[G]$ and $b \in \mathbb{C}[[G]]$ one
        has
        \begin{equation}
            \label{eq:normBimoduleStuff}
            \norm{ab}_{\ell^\infty, L, \sigma, -R}
            \le
            \norm{a}_{L, \sigma, R}
            \norm{b}_{\ell^\infty, L, \sigma, -R}
            \quad
            \textrm{as well as}
            \quad
            \norm{ba}_{\ell^\infty, L, \sigma, -R}
            \le
            \norm{a}_{L, \sigma, R}
            \norm{b}_{\ell^\infty, L, \sigma, -R}.
        \end{equation}
    \item \label{item:ellInftyBimodule} The bimodule structure of
        $\mathbb{C}[[G]]$ with respect to the group algebra
        $\mathbb{C}[G]$ restricts to a continuous bimodule structure
        on $\ellInfty_{L, \sigma, -R}(G)$ with respect to the Banach
        algebra $\ellOne_{L, \sigma, R}(G)$.
    \item \label{item:AintoEllInftyBimoduleMap} The canonical
        inclusion $\ellOne_{L, \sigma, R}(G) \subseteq \ellInfty_{L,
          \sigma, -R}(G)$ is a continuous bimodule morphism.
    \end{propositionlist}
\end{proposition}
\begin{proof}
    The first part is just Lemma~\ref{lemma:BimodContinuity}. This
    implies the second part and hence the third part follows from the
    general algebraic considerations from
    Section~\ref{sec:NotationPreliminaries}.
\end{proof}
\begin{corollary}
    \label{corollary:SubCaseBimodIntoDual}%
    Let $\sigma$ be a submultiplicative growth function. Then for all
    $R \ge 0$ the algebra $\ellOne_{L, \sigma, R}(G)$ embeds as
    bimodule into its topological dual via \eqref{eq:DualIslInfty} in
    a continuous way.
\end{corollary}
Indeed, the canonical bimodule structures on the dual correspond to
the bimodule structures on $\ellInfty_{L, \sigma, -R}(G)$ as both are
the restrictions of the bimodule structures of $\mathbb{C}[[G]]$ and
$\mathbb{C}[G]^*$ with respect to $\mathbb{C}[G]$ which then can be
extended by continuity from $\mathbb{C}[G]$ to $\ellOne_{L, \sigma,
  R}(G)$.

In a next step we pass to the projective limit $\algebra{A}_\sigma(G)$
of the algebras $\ellOne_{L, \sigma, R}(G)$ for $R \longrightarrow
\infty$. First we note that the projective system is \emph{reduced},
i.e. the projective limit $\algebra{A}_\sigma(G)$ is dense in all the
spaces $\ellOne_{L, \sigma, R}(G)$:
\begin{lemma}
    \label{lemma:ProjectiveLimitReduced}%
    Let $\sigma$ be an almost submultiplicative growth function.  Then
    the projective system of Banach spaces $\{\ellOne_{L, \sigma,
      R}(G)\}_{R \ge 0}$ is reduced.
\end{lemma}
\begin{proof}
    Indeed, the subspace $\mathbb{C}[G] \subseteq
    \algebra{A}_\sigma(G)$ is already dense in each $\ellOne_{L,
      \sigma, R}(G)$.
\end{proof}

It follows that the dual space of $\algebra{A}_\sigma(G)$ is given by
the \emph{inductive limit} of the dual spaces $\ellOne_{L, \sigma,
  R}(G)'$ as a vector space, i.e. we have
\begin{equation}
    \label{eq:DualSpaceAsigma}
    \algebra{A}_\sigma(G)'
    =
    \left(
        \projlim_{R \longrightarrow \infty}
        \ellOne_{L, \sigma, R}(G)
    \right)'
    =
    \indlim_{R \longrightarrow \infty}
    \ellOne_{L, \sigma, R}(G)'
    \cong
    \bigcup_{R \ge 0}
    \ellInfty_{L, \sigma, -R}(G)
\end{equation}
as \emph{vector spaces}, see e.g. \cite[§26,
Satz~1.6]{floret.wloka:1968a}. Note however, without more detailed
information this is just an isomorphism of vector spaces, the
corresponding locally convex inductive limit of the right hand side
can still be rather complicated. In order to determine the locally
convex topology coming from this inductive limit structure, we need
the following observation:
\begin{lemma}
    \label{lemma:ProjectiveSystemCompact}%
    Let $\sigma$ be an almost submultiplicative growth function with
    $\beta_G \earlier \sigma$. Then the projective system $\ellOne_{L,
      \sigma, R}(G)\}_{R \ge 0}$ is nuclear and hence compact.
\end{lemma}
\begin{proof}
    This is essentially the Grothendieck-Pietsch criterion: we have to
    show that for every $R \ge 0$ there is a $R' > R$ such that the
    inclusion map $\ellOne_{L, \sigma, R'}(G) \longrightarrow
    \ellOne_{L, \sigma, R}(G)$ is nuclear. This inclusion map can be
    written as
    \begin{equation*}
        a
        \; \mapsto \;
        \sum_{g \in G}
        \frac{1}{\sigma(L(g))^{R' - R}}
        \frac{\basis{e}_g}{\sigma(L(g))^R}
        \sigma(L(g))^{R'} \delta_g(a),
    \end{equation*}
    where we view $\sigma(L(g))^{R'} \delta_g\colon \ellOne_{L,
      \sigma, R'}(G) \longrightarrow \mathbb{C}$ as a continuous
    linear functional of functional norm $1$ and where
    $\frac{\basis{e}_g}{\sigma(L(g))^R} \in \ellOne_{L, \sigma, R}(G)$
    is a unit vector. The fact that ($*$) is nuclear is thus a
    consequence of the summability of the pre-factors for large enough
    $R' - R$, which is the case if $\sigma_G \earlier \sigma$ by the
    very argument used in the proof of
    Theorem~\ref{theorem:Nuclear}. Finally, we use that a nuclear map
    is always compact.
\end{proof}

For a compact and reduced projective system the algebraic
identification \eqref{eq:DualSpaceAsigma} can be turned into a
topological identification if one uses the strong topology on the
dual, see \cite[§26, Satz~2.4]{floret.wloka:1968a}:
\begin{theorem}[The dual of $\algebra{A}_\sigma(G)$]
    \label{theorem:DualOfAsigma}%
    Let $\sigma$ be an almost submultiplicative growth function with
    $\beta_G \earlier \sigma$. Then the locally convex inductive limit
    topology on the inductive system $\{\ellOne_{L, \sigma,
      R}(G)'\}_{R \ge 0}$ yields the strong topology on the dual of
    $\algebra{A}_\sigma(G)$, i.e. we have
    \begin{equation}
        \label{eq:StrongDual}
        \algebra{A}_\sigma(G)'
        =
        \indlim_{R \longrightarrow \infty}
        \ellOne_{L, \sigma, R}(G)'
        \cong
        \indlim_{R \longrightarrow \infty}
        \ellInfty_{L, \sigma, -R}(G)
    \end{equation}
    in the sense of the identification \eqref{eq:DualSpaceAsigma} with
    the strong topology on the dual.
\end{theorem}
\begin{remark}
    \label{remark:StrongDual}%
    Under the assumption $\beta_G \earlier \sigma$, which we found
    already very useful in Theorem~\ref{theorem:Nuclear},  we have on the one
    hand a rather easy description of the dual space of
    $\algebra{A}_\sigma(G)$ as vector space, namely the union of all
    the spaces $\ellInfty_{L, \sigma, -R}(G)$. Recall that the locally
    convex inductive limit topology is determined by all the seminorms
    $\halbnorm{p}$ on $\bigcup_{R \ge 0} \ellInfty_{L, \sigma, -R}(G)$
    such that the restriction of $\halbnorm{p}$ to any of the
    individual spaces $\ellInfty_{L, \sigma, -R}(G)$ is continuous. In
    general, it is quite hard to write down an explicit
    characterization of such seminorms. Nevertheless, the locally
    convex inductive limit topology has the nice universal property,
    that a linear map $\Phi\colon \bigcup_{R \ge 0} \ellInfty_{L,
      \sigma, -R}(G) \longrightarrow V$ into some other locally convex
    space $V$ is continuous iff each restriction
    $\Phi\at{\ellInfty_{L, \sigma, -R}(G)}$ is continuous. On the
    other hand, the strong topology on the dual
    $\algebra{A}_\sigma(G)$ is determined by the explicit seminorms
    $\halbnorm{p}_B(\varphi) = \sup_{a \in B}|\varphi(a)|$ where $B
    \subseteq \algebra{A}_\sigma(G)$ ranges over all bounded subsets
    of $\algebra{A}_\sigma(G)$. From the general theory we know that
    the strong topology is complete and not just sequentially complete
    (as e.g. the weak$^*$ topology). It is the combination of the two
    aspects which makes the dual of $\algebra{A}_\sigma(G)$
    manageable. Note however, that the inductive limit in
    \eqref{eq:DualSpaceAsigma} will \emph{not} be a strict inductive
    limit.
\end{remark}

While in the Banach space situation the dual of an algebra carries a
continuous bimodule structure, the bimodule structure in the general
locally convex case might no longer be continuous. We recall the
following well-known statement:
\begin{proposition}
    \label{proposition:ContinuousBimoduleDual}%
    Let $\algebra{A}$ be a locally convex algebra and equip its dual
    $\algebra{A}'$ with the strong topology. Then $\algebra{A}'$ is a
    bimodule over $\algebra{A}$ via \eqref{eq:BimoduleDual} such that
    the bimodule structure is separately continuous.
\end{proposition}
Note that in the Banach space situation the strong topology is the
usual norm topology on the dual. Thus this statement generalizes the
results from Proposition~\ref{proposition:BimoduleSubCase}. Since
$\algebra{A}_\sigma(G)$ is included continuously into $\ell^1(G)$ and
hence, via \eqref{eq:InclusionsAintoEllInfty}, into the dual
$\algebra{A}_\sigma(G)'$, we obtain the following statement:
\begin{corollary}
    \label{corollary:AsigmaIntoAsigmaDual}%
    Let $\sigma$ be an almost submultiplicative growth function with
    $\beta_G \earlier \sigma$. Then $\algebra{A}_\sigma(G)$ is
    continuously included into its strong dual
    $\algebra{A}_\sigma(G)'$ as a bimodule.
\end{corollary}

%
%

\section{The Spaces $\cNull_{L, \sigma, -R}(G)$}
\label{sec:CnullSpaces}

For fixed $R$ we have the Banach spaces $\ellInfty_{L, \sigma,
  -R}(G)$, which constitute the topological duals of $\ellOne_{L,
  \sigma, R}(G)$, and the evaluation functionals $\delta_g =
\tr(S(\basis{e}_g) \argument)$ for $g \in G$ are continuous linear
functionals, i.e. elements of the topological dual. Moreover, elements
in $\ellInfty_{L, \sigma, -R}(G)$ are certain formal series in the
basis vectors $\basis{e}_g$. Thus it is reasonable to ask whether
these vectors still form a topological basis of $\ellInfty_{L, \sigma,
  -R}(G)$. This is not true for general reasons but the failure of
being a topological basis is not too bad: we have to take into account
all spaces $\ellInfty_{L, \sigma, -R}(G)$ for $R \longrightarrow
\infty$ in the end and we shall show in this section that between any
two of them is a space $\cNull_{L, \sigma, -R}(G)$ for which the basis
vectors indeed form an unconditional Schauder basis.

As for general weighted $\ellInfty$-spaces we define the following
subset of $\ellInfty_{L, \sigma, -R}(G)$
\begin{equation}
    \label{eq:cNullDef}
    \cNull_{L, \sigma, -R}(G)
    =
    \bigg\{
    b \in \mathbb{C}[[G]]
    \; \bigg| \;
    \textrm{for all } \epsilon > 0
    \textrm{ there is a finite subset }
    K \subseteq G
    \textrm{ with }
    \sup_{g \in G \setminus K} \frac{|b_g|}{\sigma(L(g))^R}
    < \epsilon
    \bigg\},
\end{equation}
see also \eqref{eq:cNUllmu}. The idea is that this space corresponds
to the usual sequence space of zero sequences $\cNull$ equipped with
the supremum norm, viewed as a closed subspace of $\ell^\infty$. Now
the sequences do of course not converge to $0$ but require the weights
for the estimate of being small outside a finite subset. From the
general theory we obtain the following important statement:
\begin{proposition}
    \label{proposition:cNullProperties}%
    Let $R \ge 0$.
    \begin{propositionlist}
    \item \label{item:cNullClosedSubspace} With respect to the norm
        $\norm{\argument}_{\ell^\infty, L, \sigma, -R}$ the set
        $\cNull_{L, \sigma, -R}(G)$ is a closed subspace of
        $\ellInfty_{L, \sigma, -R}(G)$.
    \item The $\mathbb{C}$-span of the vectors $\{\basis{e}_g\}_{g \in
          G}$ is dense in $\cNull_{L, \sigma, -R}(G)$.
    \item The vectors $\{\basis{e}_g\}_{g \in G}$ form an
        unconditional Schauder basis in $\cNull_{L, \sigma, -R}(G)$.
    \item For $R' > R$ we have the continuous inclusion
        \begin{equation}
            \label{eq:ellInftyIncNull}
            \ellInfty_{L, \sigma, -R}(G)
            \subseteq
            \cNull_{L, \sigma, -R'}(G)
        \end{equation}
        with dense image.
    \end{propositionlist}
\end{proposition}
\begin{proof}
    The first three statements are true in general and folklore, see
    Appendix~\ref{sec:LocallyConvex}. For the last statement, we
    already know that $\ellInfty_{L, \sigma, -R}(G) \subseteq
    \ellInfty_{L, \sigma, -R'}(G)$ is a continuous inclusion. Now let
    $b \in \ellInfty_{L, \sigma, -R}(G)$ be given. Then for all $g \in
    G$ we have
    \[
    \frac{|b_g|}{\sigma(L(g))^{R'}}
    =
    \frac{|b_g|}{\sigma(L(g))^R} \frac{1}{\sigma(L(g))^{R' - R}}
    \le
    \norm{b}_{\ell^\infty, L, \sigma, -R}
    \frac{1}{\sigma(L(g))^{R' - R}}.
    \]
    Since $\sigma$ is unbounded and monotonic by definition of a
    growth function, the sequence $(\frac{1}{\sigma(n)})_{n \in
      \mathbb{N}_0}$ is a monotonic zero sequence. Moreover, for every
    $n \in \mathbb{N}_0$ there are only finitely many group elements
    $g \in G$ with $L(g) = n$, namely $\sigma_G(n)$ many. Hence for
    $\epsilon > 0$ the subset
    \[
    K_\epsilon
    =
    \left\{
        g \in G
        \; \Big| \;
        \tfrac{1}{\sigma(L(g))^{R' - R}} \ge \epsilon
    \right\}
    \]
    is finite. This shows that for $g \in G \setminus K_\epsilon$ we
    have $\frac{|b_g|}{\sigma(L(g))^{R'}} < \epsilon
    \norm{b}_{\ell^\infty, L, \sigma, -R}$. From here we can conclude
    \eqref{eq:ellInftyIncNull}. Since already $\mathbb{C}[G]$ is dense
    in $\cNull_{L, \sigma, -R'}(G)$, the density of $\ellInfty_{L,
      \sigma, -R}(G)$ follows, too.
\end{proof}

In particular, the inductive limits of the Banach spaces $\cNull_{L,
  \sigma, -R}(G)$ and of the Banach spaces $\ellInfty_{L, \sigma,
  -R}(G)$ for $R \longrightarrow \infty$ coincide. Hence it will
sometimes be advantageous to replace the $\ellInfty_{L, \sigma,
  -R}(G)$ by the easier $\cNull_{L, \sigma, -R}(G)$ where we can use
e.g. the existence of an unconditional Schauder basis.

Concerning the algebraic structure we have to restrict again to the
submultiplicative case since only here we have an algebra structure
for $\ellOne_{L, \sigma, R}(G)$. Hence it makes sense to ask whether
$\cNull_{L, \sigma, -R}(G) \subseteq \ellInfty_{L, \sigma, -R}(G)$ is
a sub bimodule. Indeed, this is the case:
\begin{theorem}[The bimodule $\cNull_{L, \sigma, -R}(G)$]
    \label{theorem:cNullSubBimodule}%
    Let $\sigma$ be a submultiplicative growth function. Then
    $\cNull_{L, \sigma, -R}(G)$ is a sub bimodule of $\ellInfty_{L,
      \sigma, -R}(G)$ over $\ellOne_{L, \sigma, R}(G)$.
\end{theorem}
\begin{proof}
    Let $b \in \cNull_{L, \sigma, -R}(G)$ and $a \in \ellOne_{L,
      \sigma, R}(G)$ be given. We have to show that $ab, ba \in
    \cNull_{L, \sigma, -R}(G)$. Let $\epsilon > 0$ and fix a finite
    subset $K \subseteq G$ with $\frac{|b_g|}{\sigma(L(g))^R} <
    \epsilon$ for $g \in G \setminus K$. Moreover, let $\tilde{K}
    \subseteq G$ be a finite subset such that
    \[
    \sum_{g \in G \setminus \tilde{K}} |a_g| \sigma(L(g))^R
    <
    \epsilon,
    \]
    which we can arrange since the series $\norm{a}_{L, \sigma, R} <
    \infty$ converges (absolutely). Finally, let $K' \subseteq G$ be
    the subset of those elements $h$ for which there exists a $g \in
    \tilde{K}$ with $g^{-1} h \in K$. Again, this subset $K'$ is
    finite.Now for  $h \in G \setminus K'$ we estimate
    \begin{align*}
        \frac{|(ab)_h|}{\sigma(L(h))^R}
        &\le
        \sum_{g \in G}
        |a_g|
        \frac{\sigma(L(g^{-1}h))^R}{\sigma(L(h))^R}
        \frac{|b_{g^{-1}h}|}{\sigma(L(g^{-1}h))^R} \\
        &\le
        \sum_{g \in G}
        |a_g| \sigma(L(g))^R
        \frac{|b_{g^{-1}h}|}{\sigma(L(g^{-1}h))^R} \\
        &\le
        \sum_{g \in G, g^{-1}h \in K}
        |a_g| \sigma(L(g))^R
        \norm{b}_{\ell^\infty, L, \sigma, -R}
        +
        \sum_{g \in G, g^{-1}h \in G \setminus K}
        |a_g| \sigma(L(g))^R
        \epsilon \\
        &\le
        \norm{b}_{\ell^\infty, L, \sigma, -R}
        \sum_{g \in G \setminus \tilde{K}}
        |a_g| \sigma(L(g))^R
        +
        \epsilon
        \norm{a}_{L, \sigma, R} \\
        &\le
        \epsilon \left(
            \norm{b}_{\ell^\infty, L, \sigma, -R}
            +
            \norm{a}_{L, \sigma, R}
        \right).
    \end{align*}
    Rescaling $\epsilon$ gives the result that outside of a finite
    subset the coefficient $\frac{|(ab)_h|}{\sigma(L(h))^R}$ is less
    than $\epsilon$. But this means $ab \in \cNull_{L, \sigma,
      -R}(G)$. The argument for $ba$ is analogous.
\end{proof}
\begin{remark}
    \label{remark:cNullAlmostSub}%
    In the almost submultiplicative case things are less easy. Here we
    first have to pass to a projective limit $\ellOne_{L, \sigma,
      R^-}(G)$ in order to have an algebra at all. Then this requires
    an inductive limit $\ellInfty_{L, \sigma, -R^-}(G)$ to get a
    bimodule structure and the dual. Inside this, one can consider the
    inductive limit $\cNull_{L, \sigma, -R^-}(G)$. Repeating the above
    estimates shows that this is indeed a sub bimodule.
\end{remark}
\begin{remark}
    \label{remark:AsigmaSubBimodules}%
    Since $\algebra{A}_\sigma(G) \subseteq \ellOne_{L, \sigma, R}(G)$
    is a subalgebra, it follows that all the subspaces $\cNull_{L,
      \sigma, R}(G)$ of the dual $\algebra{A}_\sigma(G)'$ are sub
    bimodules over $\algebra{A}_\sigma(G)$. Hence the dual of
    $\algebra{A}_\sigma(G)$ has a quite rich sub bimodule structure.
\end{remark}

%
%

\section{An Application: The Complete Growth}
\label{sec:TotalGrowth}

Recall that the \emph{complete growth function} of $G$ with respect to
a chosen word length is defined by
\begin{equation}
    \label{eq:FormalCompleteGrowth}
    F_L(z)
    = \sum_{g \in G} \basis{e}_g z^{L(g)}
    = \sum_{n = 0}^\infty
    \left(\sum_{g \in G, \sigma_G(g) = n} \basis{e}_g\right)
    z^n
    \in
    \mathbb{Z}[G][[z]],
\end{equation}
see \cite{grigorchuk.nagnibeda:1997a} where many non-trivial
properties of the complete growth were studied.

In general, this will be just a formal power series in $z$. The
purpose of this section is to find an analytic scenario where we can
make sense out of the convergence of $F_L(z)$. The first try would be
to look at our algebras $\ellOne_{L, \sigma, R}(G)$ for a suitable
choice of $\sigma$ and $R$. However, the coefficients in $F_L(z)$
consist of sums of group elements without any decaying
pre-factors. Hence we do not expect to have reasonable convergence
properties in the topologies of $\ellOne_{L, \sigma, R}(G)$, let alone
in $\algebra{A}_\sigma(G)$ where we have to require a fast decay,
controlled by $\sigma$.

Thus the next idea is to consider $F_L(z)$ as an element in the
\emph{dual} of $\algebra{A}_\sigma(G)$: here we can have a not to fast
growth of the coefficients. Again, this supports the point of view to
think of the duality between test functions $\algebra{A}_\sigma(G)$
and distributions and hence of a convergence of
\eqref{eq:FormalCompleteGrowth} in a distributional sense.
Under the nuclearity assumption this turns out to be the case:
\begin{theorem}[Convergence of the complete growth]
    \label{theorem:CompleteGrowthConverges}%
    Let $(1 + n) \earlier \sigma$ be an almost submultiplicative
    growth function. Then the following statements are equivalent:
    \begin{theoremlist}
    \item We have $\beta_G \earlier \sigma$.
    \item There exists an $R \ge 0$ such that $F_L(z)$ converges
        absolutely for $|z| \le 1$ with respect to
        $\norm{\argument}_{\ell^\infty, L, \sigma, -R}$.
    \end{theoremlist}
\end{theorem}
\begin{proof}
    Clearly it suffices to consider $|z| = 1$.  The absolute
    convergence of the series \eqref{eq:FormalCompleteGrowth} with
    respect to the norm $\norm{\argument}_{\ell^\infty, L, \sigma, -R}$
    for $|z| = 1$ means that
    \[
    \sum_{g \in G}
    \norm{\basis{e}_g z^{L(g)}}_{\ell^\infty, L, \sigma, -R}
    =
    \sum_{g \in G} \frac{1}{\sigma(L(g))^R}
    <
    \infty.
    \]
    But this is equivalent to say that $g \mapsto
    \frac{1}{\sigma(L(g))}$ is $\ell^p$-summable over $G$ for some $p
    \ge 1$, namely $p = R$. By Remark \ref{remark:NuclearityCondition}
    this is equivalent to $\sigma_G \earlier \sigma$ and hence
    $\beta_G \earlier \sigma$, by
    Lemma~\ref{lemma:VolumeSurfaceGrowthSmallerSigma}.
\end{proof}
\begin{corollary}
    \label{corollary:CompleteGrowth}%
    For an almost submultiplicative growth function $\sigma$ with
    $\beta_G \earlier \sigma$ the complete growth $F_L$ is a
    holomorphic function of $z$ on the unit disc taking values in the
    strong dual of $\algebra{A}_\sigma(G)$.
\end{corollary}

%
%

\appendix

%
%

\section{Some Locally Convex Analysis}
\label{sec:LocallyConvex}

In this short appendix we collect some basic and well-known results
from locally convex analysis. More details can be found in
e.g. \cite{jarchow:1981a, floret.wloka:1968a}.

We make extensive use of the tensor product of seminorms: if
$\halbnorm{p}$ and $\halbnorm{q}$ are seminorms on vector spaces $V$
and $W$, respectively, then $\halbnorm{p} \tensor \halbnorm{q}$
defined by
\begin{equation}
    \label{eq:TensorProductSeminorms}
    (\halbnorm{p} \tensor \halbnorm{q})(z)
    =
    \inf\Big\{
    \sum\nolimits_i \halbnorm{p}(v_i)\halbnorm{q}(w_i)
    \; \Big| \;
    z = \sum\nolimits_i v_i \tensor w_i
    \Big\}
\end{equation}
is a seminorm on the tensor product $V \tensor W$ which is a norm iff
$\halbnorm{p}$ and $\halbnorm{q}$ are norms. On factorizing tensors
one has $(\halbnorm{p} \tensor \halbnorm{q})(v \tensor w) =
\halbnorm{p}(v)\halbnorm{q}(w)$.

Very often we need seminorms build from $\ell^1$-like norms with
respect to  certain vector space bases and weighted counting
measures. For their tensor product one has the following result:
\begin{lemma}
    \label{lemma:TensorProductEllEinsSeminorms}%
    Let $V = \mathbb{C}\textrm{-}\spann\{e_i\}_{i \in I}$ and $W =
    \mathbb{C}\textrm{-}\spann\{f_j\}_{j \in J}$ be vector spaces with
    bases. Moreover, let $\mu_i \ge 0$ and $\nu_j \ge 0$ for $i \in I$
    and $j \in J$ be given and consider the $\ell^1$-seminorms
    \begin{equation}
        \label{eq:EllEinsMu}
        \norm{v}_{\ell^1(\mu)} = \sum_{i \in I} |v_i| \mu_i
        \quad
        \textrm{and}
        \quad
        \norm{w}_{\ell^1(\nu)} = \sum_{j \in J} |w_j| \nu_j,
    \end{equation}
    where $v = \sum_{i \in I} v_i e_i$ and $w = \sum_{j \in J} w_j
    f_j$. Then $\norm{\argument}_{\ell^1(\mu)} \tensor
    \norm{\argument}_{\ell^1(\nu)} = \norm{\argument}_{\ell^1(\mu
      \times \nu)}$, where
    \begin{equation}
        \label{eq:EllEinsMuNu}
        \norm{z}_{\ell^1(\mu \times \nu)}
        =
        \sum_{i \in I, j \in J} |z_{ij}| \mu_i \nu_j
    \end{equation}
    for $z = \sum_{i \in I, j \in J} z_{ij} e_i \tensor f_j \in V
    \tensor W$.
\end{lemma}
Note that for $\norm{\argument}_{\ell^1(\mu)}$ to be a norm we have to
require $\mu_i > 0$ for $i \in I$.

In all our constructions the existence of topological bases simplified
the arguments in a substantial way. Recall that a linearly independent
set $\{\basis{e}_i\}_{i \in I}$ in a locally convex space $V$ is
called an \emph{unconditional Schauder basis} if there are continuous
coefficient functionals $\delta_i \in V'$ such that for every $v \in
V$ there are only countably many $\delta_i(v)$ different from zero and
$v = \sum_{i \in I} \delta_i(v) \basis{e}_i$ converges: since we have
not chosen an ordering of $I$ the convergence is necessarily
unconditional. The Schauder basis is called \emph{absolute} if for
every continuous seminorm $\halbnorm{p}$ on $V$ there is another
continuous seminorm $\halbnorm{q}$ such that
\begin{equation}
    \label{eq:WhatIsSchauderBasis}
    \sum_{i \in I} |\delta_i(v)| \halbnorm{p}(\basis{e}_i)
    \le
    \halbnorm{q}(v).
\end{equation}
Note that an absolute Schauder basis is necessarily unconditional but
the converse is not true in general.
\begin{example}
    \label{example:AbsoluteSchauder}%
    A quite archetypal example of a Banach space with an absolute
    Schauder space is the sequence space $\ell^1$. Slightly more
    general, we consider a vector space $V$ with a basis
    $\{\basis{e}_i\}_{i \in I}$ and weights $\mu_i > 0$ for $i \in I$
    as we did that in
    Lemma~\ref{lemma:TensorProductEllEinsSeminorms}. Then the
    completion $\complete{V}$ of $V$ with respect to the weighted
    $\ell^1$-norm $\norm{\argument}_{\ell^1(\mu)}$ can be identified
    with the set $\complete{V}$ of formal series $v = \sum_{i \in I}
    v_i \basis{e}_i$ where for each $v$ at most countably many $v_i$
    are different from zero and $\norm{v}_{\ell^1(\mu)} = \sum_{i \in
      I} |v_i| \mu_i$ is finite. In more measure-theoretic terms the
    completion is just the space $\ell^1(I, \mu)$ where the index set
    $I$ is equipped with the weighted counting measure $\mu$.  For
    this example it is now easy to see that the evaluation functionals
    $\delta_i\colon v \mapsto v_i$ are continuous with respect to
    $\norm{\argument}_{\ell^1(\mu)}$. In fact, they have functional
    norm $\norm{\delta_i} = \frac{1}{\mu_i}$. Moreover, the former
    vector space basis $\{\basis{e}_i\}_{i \in I}$ of $V$ becomes now
    an absolute Schauder basis of its completion.
\end{example}

We also need a concrete description of the duals of weighted
$\ell^1$-spaces. Thus let again $V$ be a vector space with basis
$\{\basis{e}_i\}_{i \in I}$ and weights $\mu_i > 0$ which we assume to
be strictly positive to be in a Banach space situation for the
completion $\complete{V}$ with respect to the norm
$\norm{\argument}_{\ell^1(\mu)}$. If $\Phi\colon \complete{V}
\longrightarrow \mathbb{C}$ is a continuous linear functional then it
is completely determined by its values $\Phi(\basis{e}_i)$ on the
absolute Schauder basis. If $\norm{\Phi}$ is the functional norm of
$\Phi$ then $|\Phi(\basis{e}_i)| \le \norm{\Phi}
\norm{\basis{e}_i}_{\ell^1(\mu)} = \norm{\Phi} \mu_i$ and hence
$\norm{\Phi} \le \sup_{i \in I} \frac{|\Phi(\basis{e}_i)|}{\mu_i}$. In
fact, it is easy to see that we have equality here. Moreover, given a
linear functional on $V$, then it extends continuously to the
completion (necessarily in a unique way) iff this supremum is
finite. Thus we have
\begin{equation}
    \label{eq:DualCompleteV}
    \complete{V}'
    =
    \left\{
        \Phi \in \complete{V}^*
        \; \bigg| \;
        \sup_{i \in I} \frac{|\Phi(\basis{e}_i)|}{\mu_i} < \infty
    \right\},
\end{equation}
and the functional norm is this above supremum. Note that this
condition can also be understood as the usual supremum-norm condition
with respect to the \emph{unit vectors} $\frac{1}{\mu_i}
\basis{e}_i$. Nevertheless, for our application the basis vectors
$\basis{e}_i$ will play a more important role as we will have many
weighted counting measures $\mu$ simultaneously but keep the vectors
$\basis{e}_i$ in our application. Hence there is no point in rescaling
the vectors $\basis{e}_i$ to turn them into unit vectors for
\emph{one} of the measures $\mu$.

In general, the linear span of the coefficient functionals $\delta_i$
will not be dense in $\ell^\infty(I, \mu)$. Instead, its completion is
a much smaller subspace corresponding to the space of zero sequences
$\cNull$. In general, one considers
\begin{equation}
    \label{eq:cNUllmu}
    \cNull(I, \mu)
    =
    \Big\{
        \Phi \in \ell^\infty(I, \mu)
        \; \Big| \;
        \textrm{for all } \epsilon > 0
        \textrm{ there is a finite subset }
        K \subseteq I
        \textrm{ with }
        \sup_{i \in I\setminus K} \tfrac{1}{\mu_i} |\Phi(\basis{e}_i)|
        < \epsilon
    \Big\}.
\end{equation}
In general, this will be a proper subset of $\ell^\infty(I, \mu)$.
\begin{lemma}
    \label{lemma:cNullmu}%
    Assume $I$ is countable.  The set $\cNull(I, \mu) \subseteq
    \ell^\infty(I, \mu)$ is a closed subspace and the span of the
    evaluation functionals $\{\delta_i\}_{i \in I}$ is dense in
    $\cNull(I, \mu)$.
\end{lemma}
If $I$ is not countable, then we still can obtain the completion of
the span of the evaluation functionals where we have to add the
condition to the definition of $\cNull(I, \mu)$ that at most countably
many $\Phi(\basis{e}_i)$ are different from zero. In any case, in this
work we only need countable index sets $I$.

%
%


%
%

\end{document}
